\documentclass[3p,12pt]{elsarticle}

\usepackage{graphicx}
\usepackage{epstopdf}
\usepackage{subcaption}
\usepackage{caption}
\usepackage{amsmath}
\usepackage{amssymb}
\usepackage{balance}
\usepackage{color}
\usepackage[version=4]{mhchem}

\journal{Applied Energy}

\begin{document}

\begin{frontmatter}

\title{eXogenous Kalman Filter for Lithium-Ion Batteries State-of-Charge Estimation in Electric Vehicles}

\author{Agus Hasan$^\dagger$\footnote{Corresponding author at: Center for Unmanned Aerial Vehicles, M\ae rsk McKinney M\o ller Institute, University of Southern Denmark, Campusvej 55, Odense 5230, Denmark. Email: agha@mmmi.sdu.dk (A. Hasan).}, Martin Skriver$^\dagger$, and Tor Arne Johansen$^\ddagger$}

\address{$^\dagger$M\ae rsk McKinney M\o ller Institute, University of Southern Denmark, Odense, Denmark.\\$^\ddagger$Department of Cybernetics Engineering, Norwegian University of Science and Technology, Trondheim, Norway.}

\begin{abstract}
This paper presents a novel framework for state-of-charge estimation of rechargeable batteries in electric vehicles using a two-stage nonlinear estimator called the eXogenous Kalman filter (XKF). The nonlinear estimator consists of a cascade of nonlinear observer (NLO) and linearized Kalman filter (LKF). The NLO is used to produce a globally convergent auxiliary state estimate that is used to generate a linearized model in the time-varying Kalman filter algorithm. To demonstrate the proposed approach, we present a model of a lithium-ion battery from an equivalent circuit model (ECM). The model has linear process equations and a nonlinear output voltage equation. The method is tested using experimental data of a lithium iron phosphate (LiFePO$_4$) battery under dynamic stress test (DST) and federal urban driving schedule (FUDS). Effect on different ambient temperatures is also discussed. Compared with EKF and UKF, our proposed XKF achieve faster convergence rate, which can be attributed to the use of the NLO.
\end{abstract}

\begin{keyword}
Lithium-ion batteries, Electric Vehicles, State-of-Charge, Kalman filter, Nonlinear observer.
\end{keyword}

\end{frontmatter}

\section{Introduction}

The quest to find solutions to energy storage for the renewable industry and the rise of mobile robots, electric vehicles, and drones have led to further development of new generations of rechargeable battery technologies \cite{LLu}. Lithium-ion (Li-ion) batteries are among the most popular rechargeable batteries available in the market. They have high energy density and low self-discharge compared to other type of batteries at the same size. The main drawback is the risk of unexpected ignition from the heat generated by the battery. A slightly improvement in term of the energy density for commercial batteries is offered by lithium polymer (LiPo) batteries. LiPo batteries are common choice to power Unmanned Aerial Vehicles (UAVs). Acceleration and maneuvering over time can cause the power to drop from peak to net power. Thus, reduce the performance of the battery significantly. In an extreme case, battery failure due to overheating can cause loss of control of the UAVs leading to property damage or casualties. Thus, there is an incentive to estimate the remaining battery capacity and its performance in an UAV operation accurately. Other types of rechargeable batteries notable for their high energy density include lithium-sulfur (Li-S) and sodium-ion batteries (SIB). The former is still under development but may potentially replace lithium-ion cells not only because of its higher energy density, but also because its reduced cost from the use of sulfur \cite{George}.

The safety systems for rechargeable batteries can be improved using a smart battery management system (BMS). The BMS monitors the states and parameters of the batteries, thus provides instantaneous information about when a battery needs to be recharged and allows the BMS to prolong the battery life by preventing batteries from over-charging or over discharging.  An important component of BMS is an accurate estimation of state-of-charge (SoC) and state-of-health (SoH) \cite{KF2}. The SoC is defined as the available battery capacity, measured in Ampere-hours, and expressed as a percentage of its rated capacity. The SoC parameters can be seen as thermodynamic quantities enabling one to assess the potential energy of a battery. The SoH represents battery performance to store and deliver electrical energy relative to a new battery. The SoC and SoH cannot be measured directly and has to be estimated based on measurable variables such as current and terminal voltage. The estimation is usually done by employing mathematical or physical models to describe complex nonlinear dynamic processes arising from thermodynamics, electrode kinetics, and transport phenomena \cite{1RC}. The estimate can be used to detect battery anomaly such as overheating in failure prevention system.

There are two common approaches for SoC and SoH estimations of rechargeable batteries:
\begin{itemize}
\item \textbf{Kalman filter-based approach}. This is a natural choice since the model and the measurements in rechargeable batteries are subject to noises. Kalman filter was designed to suppress the noise affecting the battery systems. The linear Kalman filter is globally exponentially stable and optimal in the sense of minimum variance, as long as the model perfectly matches the real system, the noise is white and uncorrelated, and the covariances of the noise are known. However, battery models are nonlinear. The nonlinearities came from the measurement equation. To apply Kalman filter for SoC estimation, the system has to be linearized. The use of linearized model for the covariance error is known as the Extended Kalman filter (EKF) and has been used extensively in battery SoC estimation \cite{EKF1,EKF2,EKF3,EKF4,EKF6}. However, this method come at a price since the stability and optimality of the Kalman filter could not hold anymore. Some other variants include the Unscented Kalman filter (UKF) \cite{UKF1,UKF2,UKF3} combined with genetic algorithm \cite{Rui1,Rui2} and the particle filter \cite{Par1}. However, they are expensive in terms of computational time and are still approximations without theoretical guarantees of stability and optimality.
\item \textbf{Observer-based approach}. This approach offers global asymptotic or exponential stability. The design involves a tuning parameter that can be adjusted to pursue desired performance. Usually the analysis is based on Lyapunov stability theory. This approach include nonlinear observer \cite{OB1,OB3,OB4}, sliding-mode observer \cite{Slid}, backstepping PDE observer \cite{Moura,A4,A5}, and proportional-integral observer \cite{OB2}. These methods are relatively straightforward to implement and has higher computational efficiency.
\end{itemize}

In this paper, we present a new framework for state-of-charge estimation of rechargeable batteries using a two-stage nonlinear estimator, called the eXogenous Kalman filter (XKF), developed in \cite{Tor}. The XKF is a cascade of a nonlinear observer (NLO) and a linearized Kalman filter (LKF). The aim is to utilize each strength from both nonlinear state estimation methods, such as the global stability property of the NLO and the optimality of the LKF. The main idea is to use the exogenous state estimation obtained from the NLO to generate a linearized model for the Kalman filter.

\begin{figure}[h!]
  \centering
      \includegraphics[width=1\textwidth]{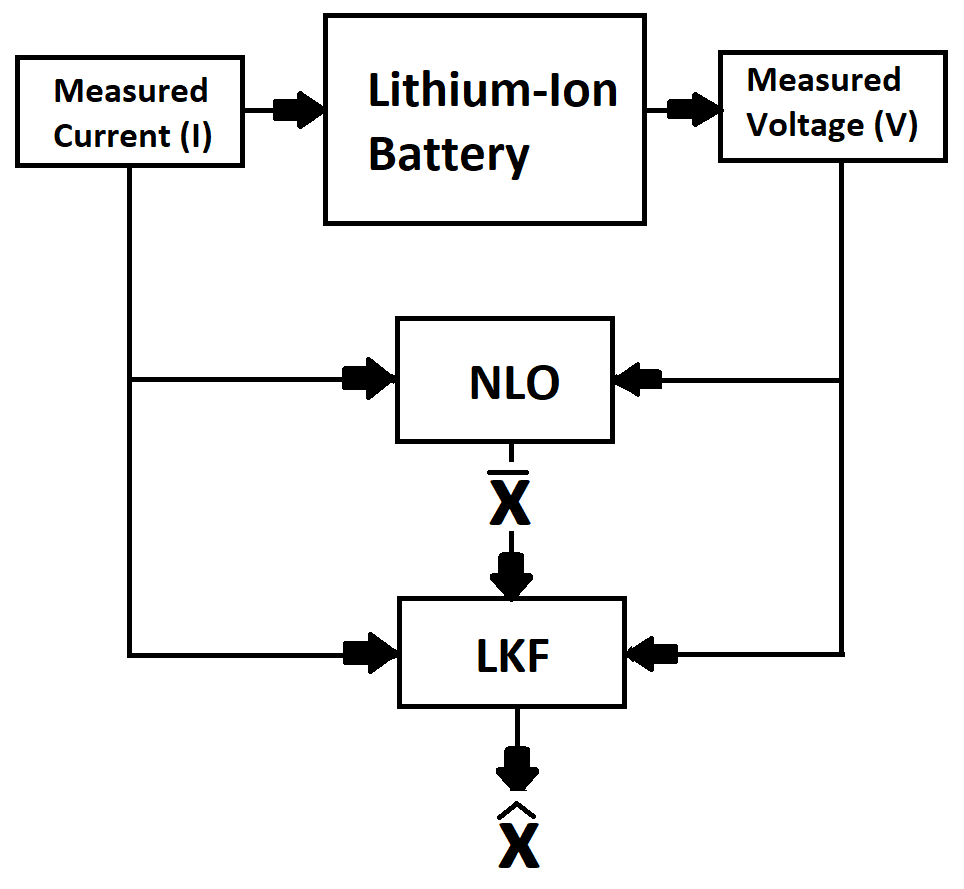}
  \caption{The schematic diagram of the XKF for SoC estimation. The XKF is a cascade of NLO and LKF.}
\label{3}
\end{figure}

In Fig. \ref{3}, we denote $\overline{\mathbf{x}}$ as the exogenous state estimation of SoC from the NLO, while $\hat{\mathbf{x}}$ denotes the SoC estimation from the LKF. The nonlinear observer gain $\mathbf{K}$ could be obtained from an analysis of the error system using Lyapunov stability theory, while the Kalman gain $\mathbf{K}_f$ is obtained after solving a Riccati equation. It has been shown that the two-stage nonlinear estimator inherits the global stability property of the NLO \cite{Tor}.

To illustrate our proposed method, we choose to work with lithium iron phosphate (LiFePO$_4$) battery. LiFePO$_4$ offers a longer life cycle and has a more stable discharge voltage than other lithium-ion batteries. It is mainly used in electric vehicles and solar energy installations \cite{Rahul}. In section 2, we present a model for LiFePO$_4$ based on Coulomb counting and 2-RC equivalent circuit model (ECM). This model is chosen because of its simple structure, thus suitable for real-time application. In section 3, we design a NLO and a LKF for our 2-RC ECM. We prove the NLO to be globally exponentially stable. Evaluation using experimental data from Dynamic Stress Test (DST) and Federal Urban Driving Schedule (FUDS) are presented in section 4. Finally, section 5 contains conclusions and future works.

\section{Lithium-ion Battery Modeling}

A lithium-ion battery uses a metal oxide cathode, a porous carbon anode, a separator, and electrolyte as conductor. The separator allows the transport of lithium-ions from one side to the other, and at the same time prevents electrons from passing through. The electrolyte is electrically conductive and contains free ions in which the lithium-ions can be transported easily \cite{John}. During discharge, the ions flow from the anode to the cathode through the electrolyte and a separator. Charging will reverse the direction and the ions flow from the cathode to the anode.

There are several types of lithium-ion batteries based on the chemistry, performance, cost, and safety. A LiFePO$_4$ battery, also called LFP battery, uses LiFePO$_4$ as a cathode material, and a graphitic carbon electrode with a metallic current collector grid as the anode. The chemical reactions in the cathode and the anode of the LiFePO$_4$ battery could be described as follow
\begin{eqnarray}
\ce{LiFePO_4 &<=>[charge][discharge]& FePO_4 + Li^+ + e^-}\nonumber\\
\ce{Li^+ + e^- + C_6 &<=>[charge][discharge]& LiC_6}\nonumber
\end{eqnarray}
The original mineral of LiFePO$_4$ is olivine which exists as a solid solution of the pure magnesium phase, forsterite, and the purely ferrous phase fayalite. So far, LiFePO$_4$ is the safest of the mainstream lithium-ion battery types. A LiFePO$_4$ battery does not require to be fully charged. Its performance even slightly improves in case of partial charge instead of a full charge. This becomes a major advantage of LiFePO$_4$ compared to lead-acid battery types. Other advantages include wide operating temperature range, excellent cycling performance, low internal resistance, and high efficiency. These are some reasons why LiFePO$_4$ has been used in electric vehicle and backup power.

A LiFePO$_4$ battery has also some drawbacks. The cell will fail if the voltage over the cell falls to less than 2.5V or increases to more than 4.2V. Furthermore, the cells do not auto-balance at the end of the charge cycle. A BMS could help to balance these cells. Furthermore, it could prevent cell under voltage by timely disconnecting the load and cell overvoltage by reducing charge current or stopping the charge process. A key component in BMS is a good model to describe the behavior of the battery under some operational condition. Lithium-ion battery modeling could be done in three different approaches:
\begin{itemize}
\item \textbf{Distributed Parameter Model}. This approach is based on thermal-electrochemical model \cite{Mou1,Hasan,Shux}. The model is derived from first principles and is written as a system of coupled partial differential equations (PDEs) and ordinary differential equations (ODEs) to describe the physics of the battery. Examples of this approach are the Doyle-Fuller-Newman model \cite{Doyle} and the single particle model (SPM) \cite{Wh}. This approach can accurately match experimental data but are complex and time consuming to solve. Advanced battery modeling using coupled PDEs-ODEs has received much attention in recent years since the estimation method for this kind of systems has been well established \cite{A1,A2,A20,A30}.
\item \textbf{Lumped Model}. This approach considers a battery as an equivalent circuit that represents the electrical characteristics of the cell \cite{Eval}. The model, which is taken from the RC ECMs, serves as a proxy model and is written as a system of ODEs. The lumped model can be solved quickly and used in real-time environments \cite{A3}. The main drawback is that this approach does not accurately capture the physics of the battery at high discharge conditions.
\item \textbf{Data-driven Model}. The last approach is based on data-driven and relies on machine learning algorithms to develop relationships between sensor inputs and output of the battery \cite{Neu1}. This approach usually requires significant amounts of hand-crafted structured training data and do not learn incrementally or interactively in real time. Next generation of machine learning technology is being researched.
\end{itemize}

\begin{figure}[h!]
  \centering
      \includegraphics[width=1\textwidth]{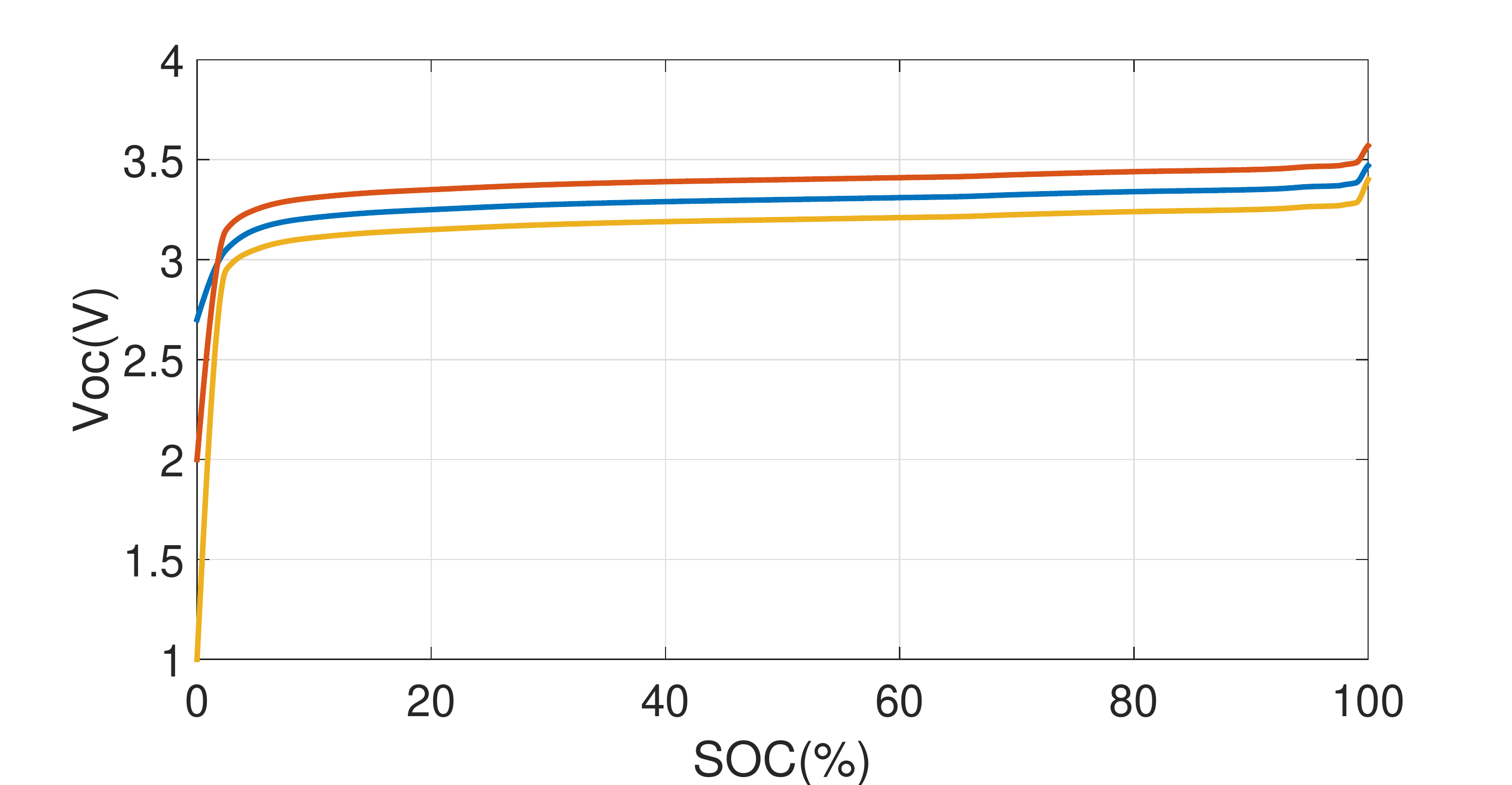}
  \caption{Typical OCV-SoC lookup tables for LiFePO$_4$ battery. The ambient temperature alters the OCV-SoC lookup tables.}
\label{1}
\end{figure}

In mobile platforms like unmanned aerial vehicles (UAVs), the complexity and accuracy of the model should be well balanced \cite{AThor}. Therefore, we present a lithium-ion battery model based on the ECM. An ECM uses a circuit consisting of voltage sources, resistors, and capacitors to simulate the dynamic characteristics of batteries \cite{Zhong}. The ECM offers a low complexity description of the lithium-ion battery dynamics with fewer states and parameters, thus more applicable in terms of practicality compared to the high-fidelity distributed parameter model. The battery SoC quantifies the usable energy at the present cycle and can be calculated by
\begin{eqnarray}
\text{SoC}(t) &=& \text{SoC}(t_0) - \frac{\int_{t_0}^{t} I(\tau)\mathrm{d}\tau}{Q_c}\label{main}
\end{eqnarray}
where $t_0$ denotes the initial time, $I$ is the electrical current, and $Q_c$ is the nominal capacity. Calculating the SoC from the measured discharging current and integrating it over time is known as Coulomb counting. This method, however, suffers from long-term drift and lack of a reference point. Therefore, the SoC should be calibrated on regular basis. Measured terminal voltages can be used for calibration. The terminal voltage can be obtained from the open circuit voltage (OCV), which underlying physical phenomenon of lithium-ion intercalation/deintercalation process and is expressed as a nonlinear function of the SoC. In practice, we create an OCV-SoC curve from experiments and use it either as an analytical expression or as a lookup table for the model. Expressing the OCV-SoC curve analytically using polynomial function is prone to incorrect estimation, since this approach is very sensitive with respect to the polynomial coefficient. Therefore, we chose to use the OCV-SoC as a lookup table. To calculate its derivative for the LKF, we approximate it using finite difference method from its interpolation. An example of an OCV-SoC curve for LiFePO$_4$ batteries could be seen in Fig. \ref{1}. The OCV-SoC curve is a function of the ambient temperature. Note that a small inaccuracy in the OCV-SoC curve may cause a large deviation in the SoC estimation, since there is a wide flat region on the OCV-SoC curve \cite{Weng}. Thus, improving the OCV models can contribute significantly to increase the accuracy of SoC estimation.

\begin{figure}[h!]
  \centering
      \includegraphics[width=1\textwidth]{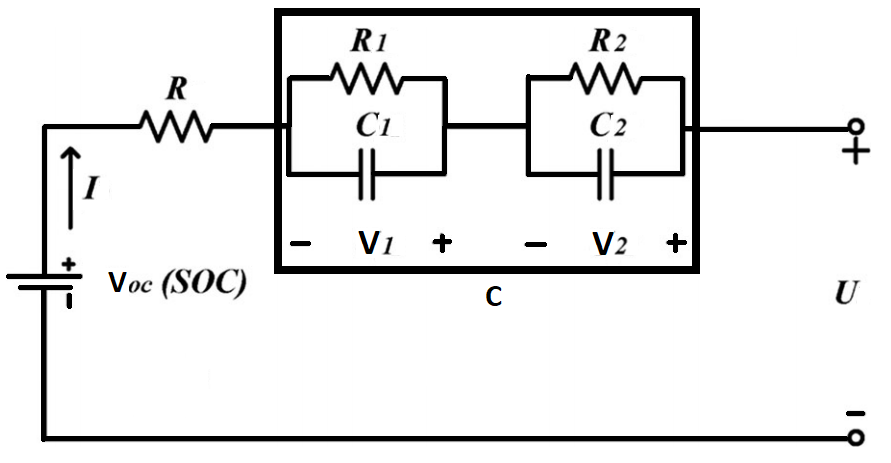}
  \caption{Schematic diagram of the 2-RC circuit for LiFePO$_4$ battery.}
  \label{2}
\end{figure}

Using the OCV-SoC lookup table and the ECM (Fig. \ref{2}), the terminal voltage $V$ can be calculated using the following formula
\begin{eqnarray}
V(t) &=& V_{oc}(\text{SoC})-I(t)R-C(t)\label{main2}
\end{eqnarray}
where $R$ is the internal resistance or the ohmic resistance. This ohmic resistance is used to represent the electrical resistance of battery components with the accumulation and dissipation of charge in the electrical double-layer. It could be easily observed from Fig. \ref{1} that $\frac{\mathrm{d} V_{oc}}{\mathrm{d}\text{SoC}}>0$. Here, $C$ is a correction factor due to model inaccuracy and environmental conditions, e.g., ambient temperature variations. Some authors use a resistor model and consider the value of $C$ as a function of ambient temperature \cite{MP1}, which could be determined using least-square fitting from experimental data. Other authors consider $C$ as a 1-RC circuit model \cite{Xhu} and 2-RC circuit model \cite{XChen}. It could be also determined using a combination of Thevenin-based ECM \cite{tuj} with the hysteresis voltage dynamics \cite{EKF5}, which offers a grasp of dynamic current-voltage characteristics and compensates the static current-voltage property.

To balance between the computational effort and accuracy, in this paper we use the 2-RC circuit model, i.e., $C(t)=V_1(t)+V_2(t)$, where
\begin{eqnarray}
\frac{\mathrm{d}V_1}{\mathrm{d}t}(t) &=& -\frac{V_1(t)}{R_1C_1}+\frac{I(t)}{C_1}\label{main3}\\
\frac{\mathrm{d}V_2}{\mathrm{d}t}(t) &=& -\frac{V_2(t)}{R_2C_2}+\frac{I(t)}{C_2}\label{main4}
\end{eqnarray}
where $R_1,R_2$ and $C_1,C_2$ are diffusion resistance and diffusion capacitance for the RC network, respectively. These parameters could also be interpreted as the mass transport effects and dynamic voltage performance. These parameters together with the ohmic resistance could be determined from the exponential-function fitting method or a simple least-square algorithm \cite{Waag}.

\section{SoC Estimation using eXogenous Kalman Filter}

Given the nonlinear battery model and measurement equations \eqref{main}-\eqref{main4}, the SoC estimation can be formulated as a nonlinear state estimation problem.

\subsection{State space model}

Differentiating \eqref{main} with respect to $t$, we have
\begin{eqnarray}
\frac{\mathrm{d}\text{SoC}}{\mathrm{d}t}(t) &=& -\frac{I(t)}{Q_c}\label{main1}
\end{eqnarray}
From \eqref{main2}-\eqref{main1}, the state-space model for the lithium-ion batteries is given by
\begin{eqnarray}
\frac{\mathrm{d}}{\mathrm{d}t}\begin{pmatrix}V_1(t)\\V_2(t)\\\text{SoC}(t)\end{pmatrix}&=&
\begin{pmatrix}-\frac{1}{R_1C_1}&0&0\\0&-\frac{1}{R_2C_2}&0\\0&0&0\end{pmatrix}\begin{pmatrix}V_1(t)\\V_2(t)\\\text{SoC}(t)\end{pmatrix}+\begin{pmatrix}\frac{1}{C_1}\\\frac{1}{C_2}\\-\frac{1}{Q_c}\end{pmatrix}I(t)\label{ss1}\\
V(t) &=& V_{oc}(\text{SoC})-V_1(t)-V_2(t)-RI(t)\label{ss2}
\end{eqnarray}
To simplify the presentation, we write \eqref{ss1}-\eqref{ss2} as
\begin{eqnarray}
\frac{\mathrm{d}\mathbf{x}}{\mathrm{d}t}(t) &=& \mathbf{A}\mathbf{x}(t)+\mathbf{B}I(t)\label{mi1}\\
V(t) &=& h(\mathbf{x})-RI(t)\label{mi2}
\end{eqnarray}
where
\begin{eqnarray}
\mathbf{x} &=& \begin{pmatrix}V_1(t)\\V_2(t)\\\text{SoC}(t)\end{pmatrix}\\
\mathbf{A} &=& \begin{pmatrix}-\frac{1}{R_1C_1}&0&0\\0&-\frac{1}{R_2C_2}&0\\0&0&0\end{pmatrix}\\
\mathbf{B} &=& \begin{pmatrix}\frac{1}{C_1}\\\frac{1}{C_2}\\-\frac{1}{Q_c}\end{pmatrix}\\
h(\mathbf{x}) &=& V_{oc}(\text{SoC})-V_1(t)-V_2(t)
\end{eqnarray}
To incorporate the model and measurement uncertainties and inaccuracies, noises are added into \eqref{mi1}-\eqref{mi2}, thus the complete model becomes
\begin{eqnarray}
\frac{\mathrm{d}\mathbf{x}}{\mathrm{d}t}(t) &=& \mathbf{A}\mathbf{x}(t)+\mathbf{B}I(t)+w(t)\\
V(t) &=& h(\mathbf{x})-RI(t)+v(t)\\
w &\sim& (0,\mathbf{Q}_n)\\
v &\sim& (0,\mathbf{R}_n)
\end{eqnarray}
where $w$ denotes the process noise, $\mathbf{Q}_n$ is the process noise covariance, $v$ is the measurement noise, and $\mathbf{R}_n$ is the measurement noise covariance. To use the optimal filter, the covariances of the noise $\mathbf{Q}_n$ and $\mathbf{R}_n$ have to be known. These values can be obtained using Bayesian, maximum likelihood, covariance matching, and correlation techniques \cite{Row}. Furthermore, the linear Kalman filter is optimal under the assumption that $w$ and $v$ are white noises and uncorrelated. In this section, we design the SoC estimation using the XKF, which is a cascade of a NLO and LKF.

\subsection{Nonlinear Observer (NLO)}

The NLO is designed as follow
\begin{eqnarray}
\frac{\mathrm{d}\overline{\mathbf{x}}}{\mathrm{d}t}(t)(t) &=& \mathbf{A}\overline{\mathbf{x}}(t)+\mathbf{B}I(t)+\mathbf{K}\left(V(t)-\overline{V}(t)\right)\label{no1}\\
\overline{V}(t) &=& h(\overline{\mathbf{x}})-RI(t)\label{no2}
\end{eqnarray}
where $\overline{\mathbf{x}}$ denotes the exogenous state estimation from the NLO and $\mathbf{K}=\begin{pmatrix}k_1&k_2&k_3\end{pmatrix}^{\intercal}\in\mathbb{R}^3$ is the observer gain to be determined later. Let us define the difference between the actual state and the measured state as $\breve{\mathbf{x}}(t) = \mathbf{x}(t) - \overline{\mathbf{x}}(t)$ and the difference between the actual terminal voltage and the measured terminal voltage as $\breve{V}(t) = V(t)-\overline{V}(t)$. The error equation is obtained as follow
\begin{eqnarray}
\frac{\mathrm{d}\breve{\mathbf{x}}}{\mathrm{d}t}(t) &=& \mathbf{A}\breve{\mathbf{x}}(t)-\mathbf{K}\breve{V}(t)\label{err}\\
\breve{V}(t) &=& h(\mathbf{x}) - h(\overline{\mathbf{x}})\label{meas}
\end{eqnarray}
Utilizing the smoothness property of the OCV-SoC lookup table, we can linearize the nonlinear measurement equation \eqref{meas} using the mean value theorem. Note that the continuity of the OCV-SoC lookup table can be obtained by applying interpolation from the experimented data. First, we write
\begin{eqnarray}
\breve{V}(t) &=& -\breve{V}_1(t)-\breve{V}_2(t)+V_{oc}(\text{SoC})-V_{oc}(\overline{\text{SoC}})
\end{eqnarray}
From the mean value theorem, there exist $\text{SoC}\leq\xi\leq\overline{\text{SoC}}$, such that
\begin{eqnarray}
V_{oc}(\text{SoC})-V_{oc}(\overline{\text{SoC}}) = \frac{\mathrm{d}V_{oc}}{\mathrm{d}\text{SoC}}(\xi)\breve{\text{SoC}}(t)
\end{eqnarray}
where $\breve{\text{SoC}}(t)=\text{SoC}(t)-\overline{\text{SoC}}(t)$. Let us denote
\begin{eqnarray}
\mathbf{T} &=& \begin{pmatrix}-1&-1&\frac{\mathrm{d}V_{oc}}{\mathrm{d}\text{SoC}}(\xi)\end{pmatrix}\label{TT}
\end{eqnarray}
Thus, we can write
\begin{eqnarray}
\breve{V}(t) = \mathbf{T}\breve{\mathbf{x}}(t)
\end{eqnarray}
and the error equation becomes
\begin{eqnarray}
\frac{\mathrm{d}\breve{\mathbf{x}}}{\mathrm{d}t}(t) &=& \left(\mathbf{A}-\mathbf{K}\mathbf{T}\right)\breve{\mathbf{x}}(t)
\end{eqnarray}
Note that $\frac{\mathrm{d}V_{oc}}{\mathrm{d}\text{SoC}}(\xi)>0$ and $\mathbf{T}$ is time-varying and bounded. In this particular case, let $k_1=k_2=0$. Then, $\breve{V}_1(t)=\breve{V}_1(0)e^{-\frac{1}{R_1C_1}}$ and $\breve{V}_2(t)=\breve{V}_2(0)e^{-\frac{1}{R_2C_2}}$, i.e., $\breve{V}_1(t)$ and $\breve{V}_2(t)$ are exponentially stable. Correspondingly, we have $\dot{\breve{\text{SoC}}}(t)=-k_3\frac{\mathrm{d}V_{oc}}{\mathrm{d}\text{SoC}}(\xi)\breve{\text{SoC}}(t)$. Thus, if $k_3>0$ then the nonlinear observer is globally exponentially stable, i.e., the error will decay to zero and the estimate $\overline{\mathbf{x}}$ will converge to the actual value exponentially.

\subsection{Linearized Kalman Filter (LKF)}

Now we have $\overline{\mathbf{x}}$ as an estimate of $\mathbf{x}$, which is a bounded signal given by the NLO \eqref{no1}-\eqref{no2}. We use this signal as a linearization point for the linear Kalman filter. A first-order Taylor series expansion about the trajectory $\overline{\mathbf{x}}$ gives the linearized model
\begin{eqnarray}
\frac{\mathrm{d}\mathbf{x}}{\mathrm{d}t}(t) &=& \mathbf{A}\mathbf{x}(t)+\mathbf{B}I(t)+w(t)\\
V(t) &=& h(\overline{\mathbf{x}})+\mathbf{H}(\overline{\mathbf{x}})\breve{\mathbf{x}}+r(\mathbf{x},\overline{\mathbf{x}})-RI(t)+v(t)
\end{eqnarray}
where
\begin{eqnarray}
\mathbf{H}(\overline{\mathbf{x}}) &=& \frac{\partial h}{\partial\mathbf{x}}(\overline{\mathbf{x}}) = \begin{pmatrix}-1&-1&\frac{\mathrm{d} V_{oc}}{\mathrm{d}\text{SoC}}(\overline{\text{SoC}})\end{pmatrix}
\end{eqnarray}
Remark that, since $\overline{\mathbf{x}}$ is bounded and converges to $\mathbf{x}$, we can neglect the higher-order term $r(\mathbf{x},\overline{\mathbf{x}})$ in the design since it has no consequences for stability \cite{Tor}. Thus, we can design the second stage estimator $\hat{\mathbf{x}}$ using the LKF as follows
\begin{eqnarray}
\frac{\mathrm{d}\hat{\mathbf{x}}}{\mathrm{d}t}(t) &=& \mathbf{A}\hat{\mathbf{x}}(t)+\mathbf{B}I(t)+\mathbf{K}_f(t)\left(V(t)-h(\overline{\mathbf{x}})-H(\overline{\mathbf{x}})(\hat{\mathbf{x}}-\overline{\mathbf{x}})+RI(t)\right)\label{xkf}
\end{eqnarray}
The time-varying gain $\mathbf{K}_f(t)$ is given by
\begin{eqnarray}
\mathbf{K}_f(t) = \mathbf{P}(t)\mathbf{H}^{\intercal}(\overline{\mathbf{x}})\mathbf{R}_n^{-1}
\end{eqnarray}
where the error covariance matrix $\mathbf{P}(t)$ is the solution to the Riccati equation
\begin{eqnarray}
\frac{\mathrm{d}\mathbf{P}}{\mathrm{d}t}(t) &=& \mathbf{A}\mathbf{P}(t)+\mathbf{P}(t)\mathbf{A}^{\intercal}+\mathbf{Q}_n-\mathbf{K}_f(t)\mathbf{R}_n\mathbf{K}_f^{\intercal}(t)
\end{eqnarray}
with $\mathbf{P}(0)$ symmetric and positive definite. Note that, unlike EKF, the XKF uses linear time-varying measurement model that is independent of the estimate $\hat{\mathbf{x}}$ in \eqref{xkf}. To implement the algorithm with a discrete-time Kalman filter, the model has to be discretized, for example using the Euler discretization method.

\section{Evaluation using Experimental Data}

The LiFePO$_4$ battery is tested in two dynamic loading condition tests: the dynamic stress test (DST) and the federal urban driving schedule (FUDS). The DST was used to identify the model parameters, while the FUDS was used to validate the performance of the SoC estimation. The battery specification is given in Table 1. A complete description regarding the tests is given in \cite{MP1}.
\begin{table}[h!]
\begin{center}
  \begin{tabular}{ l | c | r | r}
    \hline\hline
    Type & Nominal & Nominal & Upper and lower \\
    &voltage & capacity & cut-off voltage\\ \hline
    LiFePO$_4$ & 3.3V & 2.23Ah & 3.6V and 2.0V\\
    \hline
  \end{tabular}
  \caption{Specification of the LiFePO$_4$ battery used in the experiment.}
\end{center}
\end{table}

\subsection{Dynamic Stress Test}

\begin{figure}[h!]
\centering
\begin{subfigure}[b]{1\textwidth}
   \includegraphics[width=1\linewidth]{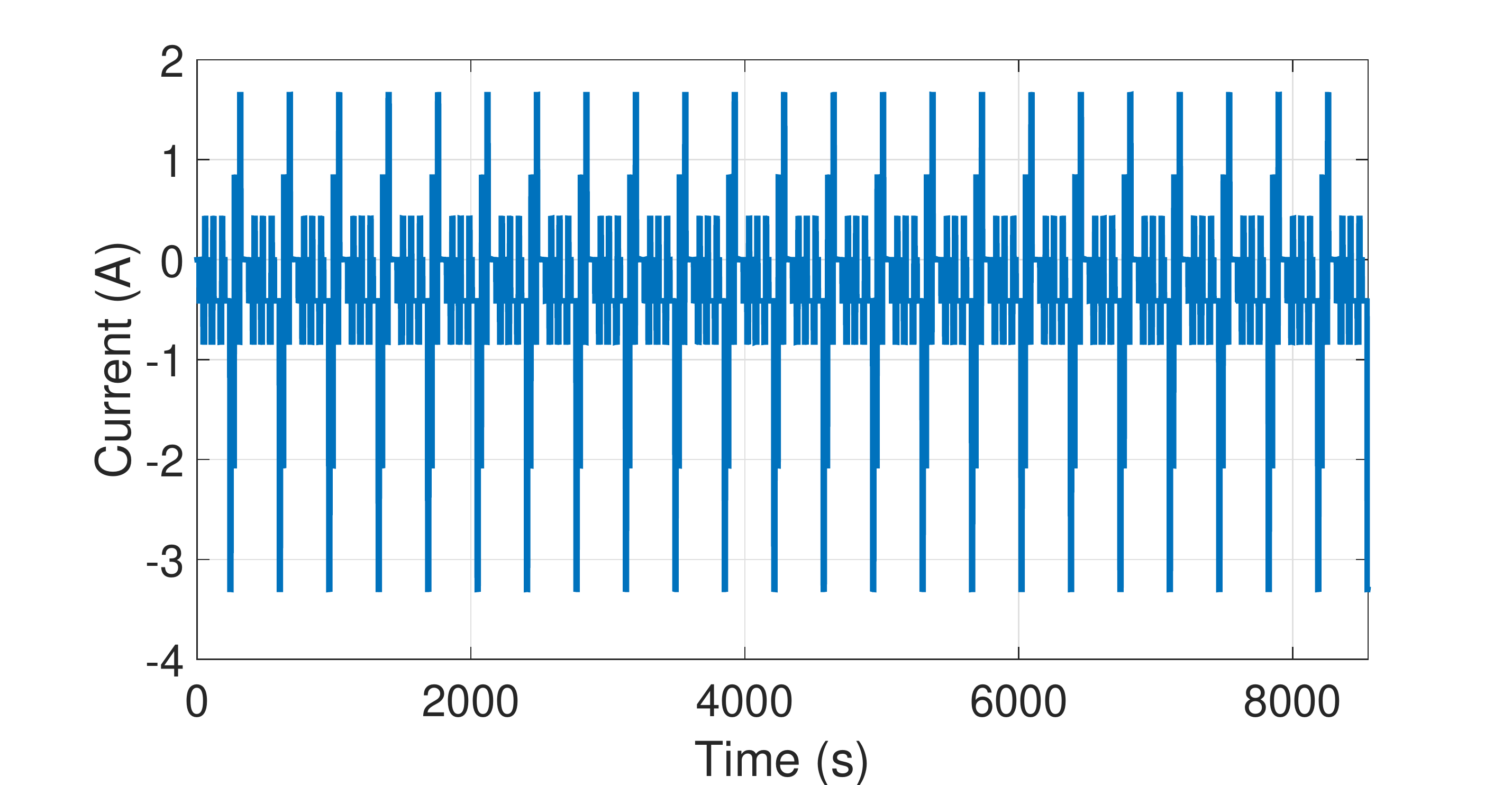}
\end{subfigure}
\begin{subfigure}[b]{1\textwidth}
   \includegraphics[width=1\linewidth]{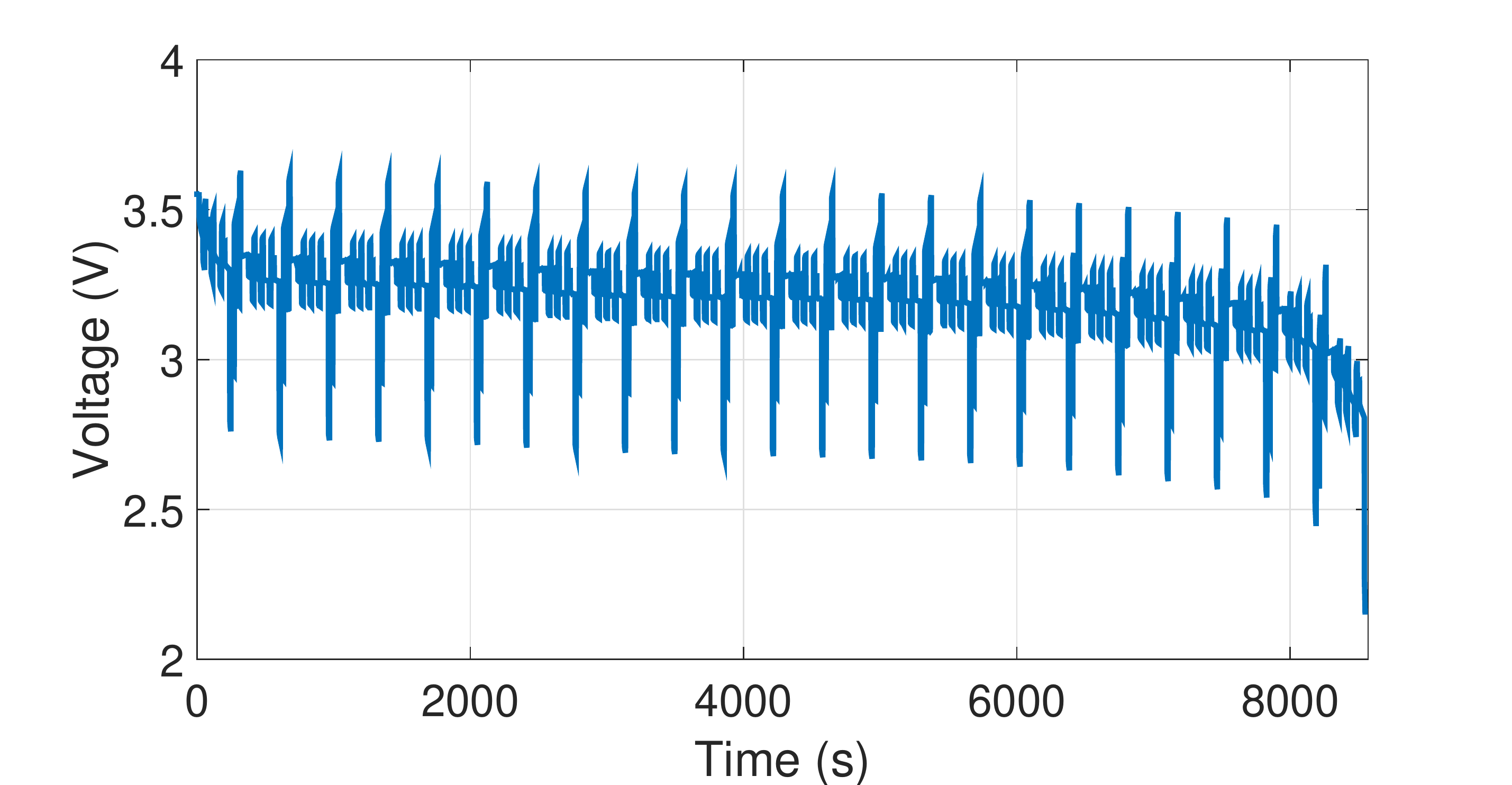}
\end{subfigure}
\caption{Measured current and terminal voltage under the dynamic stress test (DST).}
\label{4}
\end{figure}

To evaluate the validity and to identify the parameters of the battery model, the DST is performed at 20$^\circ$C. This test is designed by US Advanced Battery Consortium (USABC) to simulate a variable-power discharge regime that represents the expected demands of an electric vehicle (EV) battery. The voltage and current are measured and recorded from fully charged to empty with a sampling period of 1.0s based on the battery test bench. The accumulative charge was run continuously from 100\% SoC at 3.6V to empty at 2.0V over several cycles in a discharge process. The measured current $I$ and voltage $V$ are given in Fig. \ref{4}. This measurements are used to calibrate the estimation from the Coulomb counting method.

The root mean square (RMS) error is used to evaluate the validity of the model and the parameters and is given by
\begin{eqnarray}
\text{RMS} = \sqrt{\frac{1}{n}\Sigma_{i=1}^n\left(V_{\text{measured}}-V_{\text{estimated}}\right)^2}
\end{eqnarray}

The model parameters and the RMS error are given in Table 2. Notice that the diffusion resistance and the diffusion capacitance for $V_1$ and $V_2$ are the same. The value of $C_1$ and $C_2$ are dominant, which means the 2-RC model gave a minor effect to SoC estimation. The accuracy of the estimation is mostly determined by the OCV-SoC look up table, which its accurate value is hard to get. This becomes a main drawback when using the ECM despite of its simplicity.
\begin{table}[h!]
\begin{center}
  \begin{tabular}{ l | c | c | c | c | c}
    \hline\hline
    $R$ & $R_1$ & $C_1$ & $R_2$ & $C_2$ & \text{RMS}\\\hline
    0.18 & 0.035 & 1e6 & 0.035 & 1e6 & 8.4157e-05\\
    \hline
  \end{tabular}
  \caption{Model parameters and the RMS.}
\end{center}
\end{table}

\begin{figure}[h!]
\centering
\begin{subfigure}[b]{1\textwidth}
   \includegraphics[width=1\linewidth]{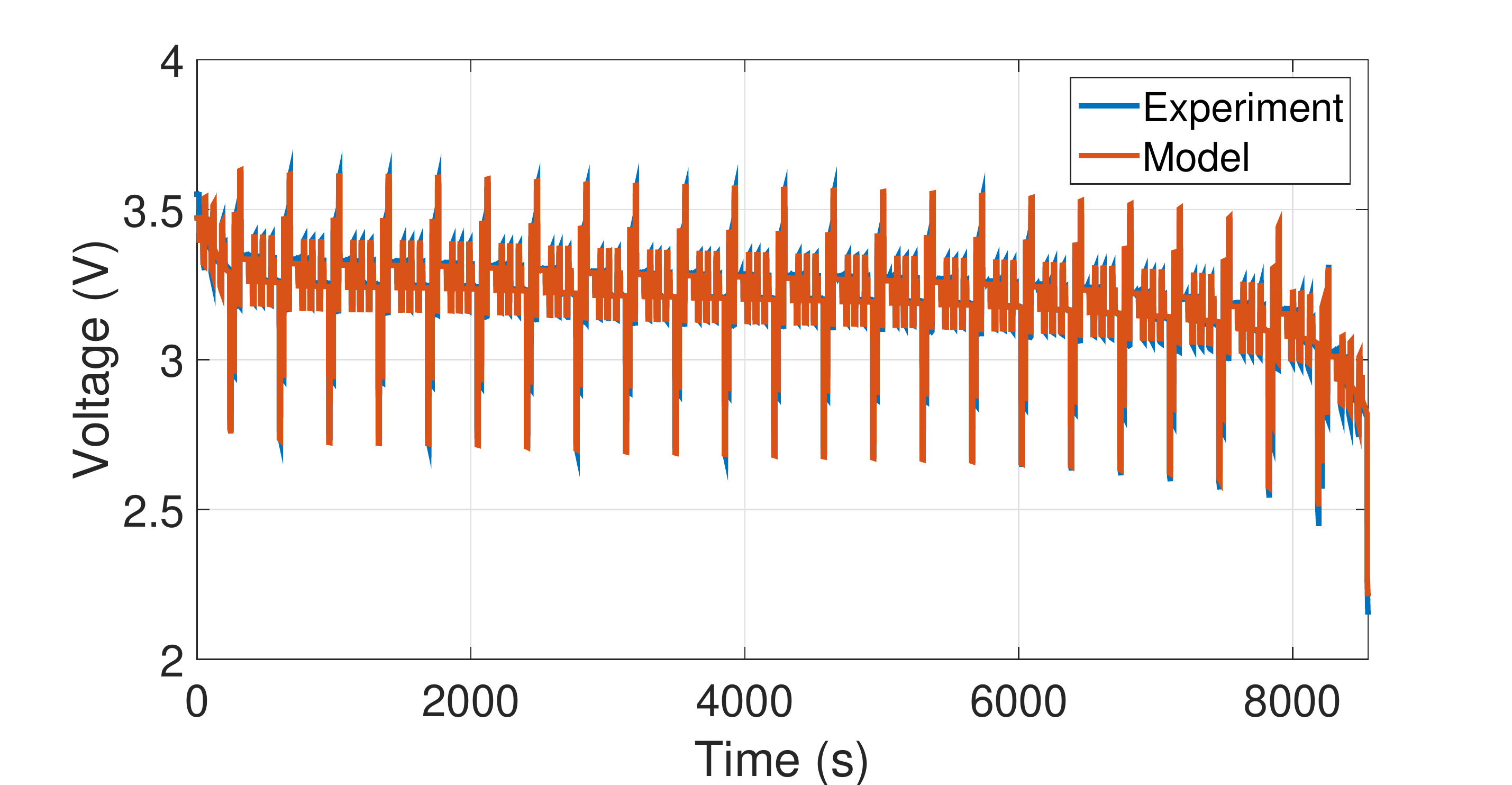}
\end{subfigure}
\begin{subfigure}[b]{1\textwidth}
   \includegraphics[width=1\linewidth]{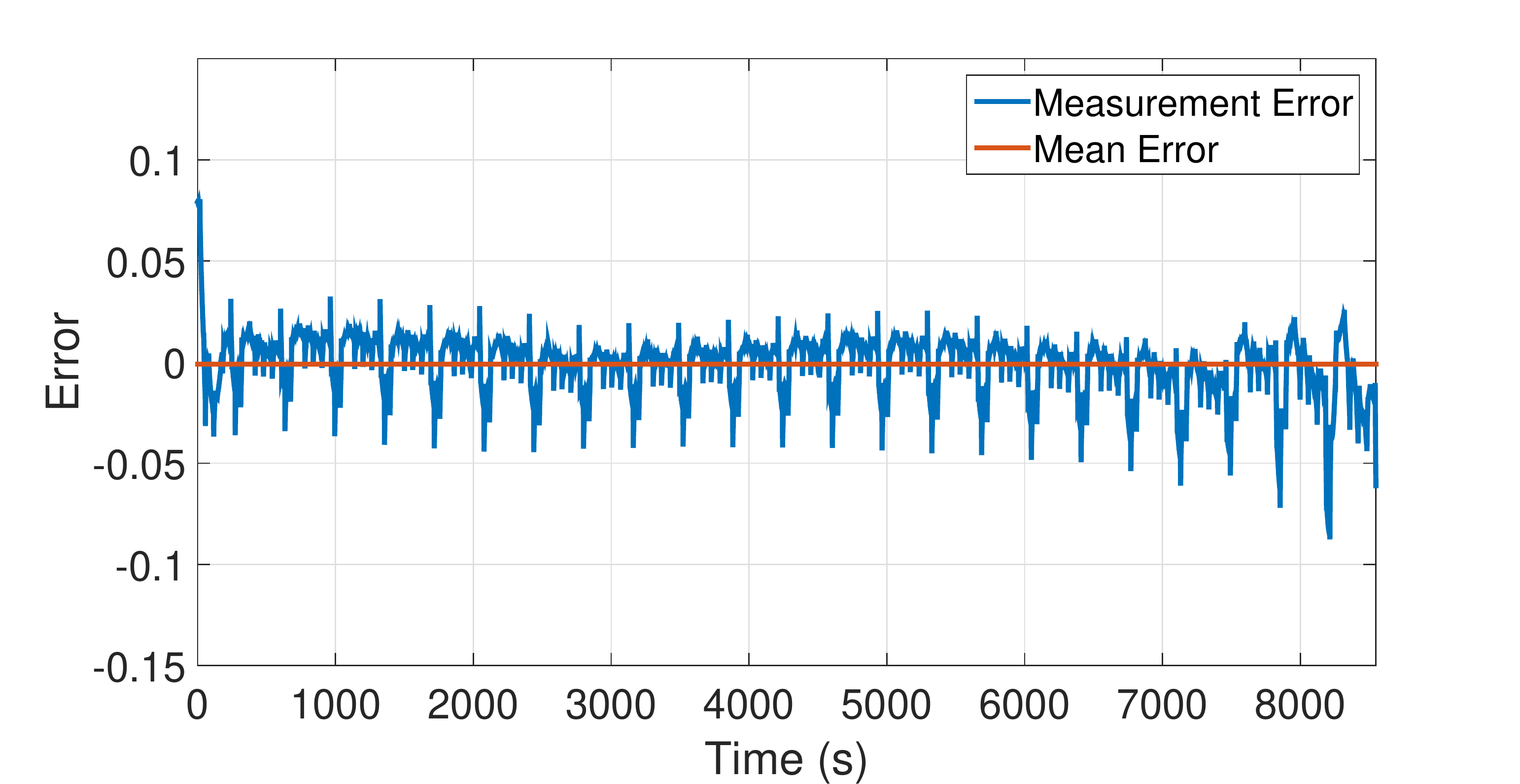}
\end{subfigure}
\caption{Measured and estimated voltage response under the DST and the model error.}
\label{5}
\end{figure}

\begin{figure}[h!]
\centering
\begin{subfigure}[b]{1\textwidth}
   \includegraphics[width=1\linewidth]{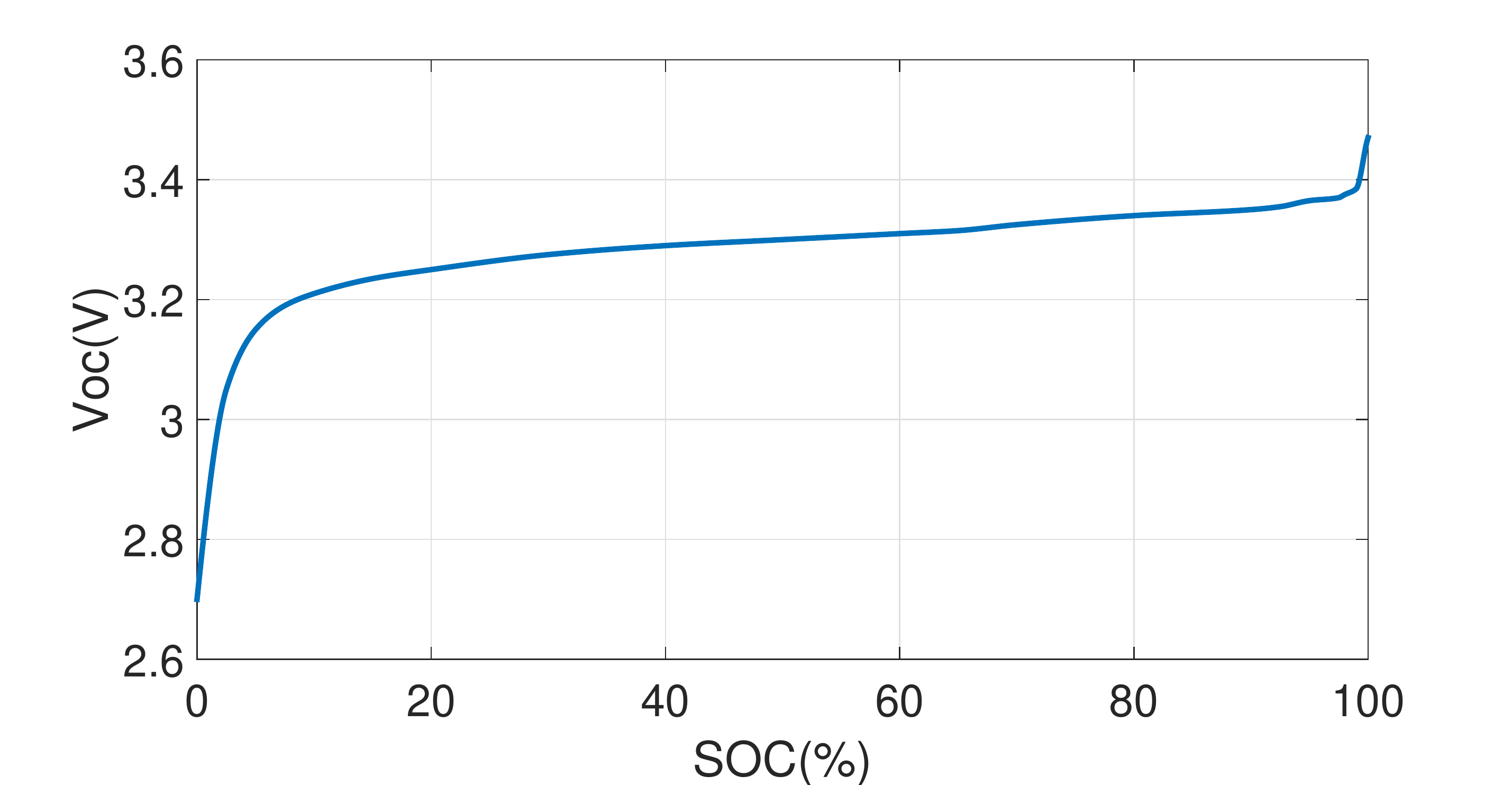}
\end{subfigure}
\begin{subfigure}[b]{1\textwidth}
   \includegraphics[width=1\linewidth]{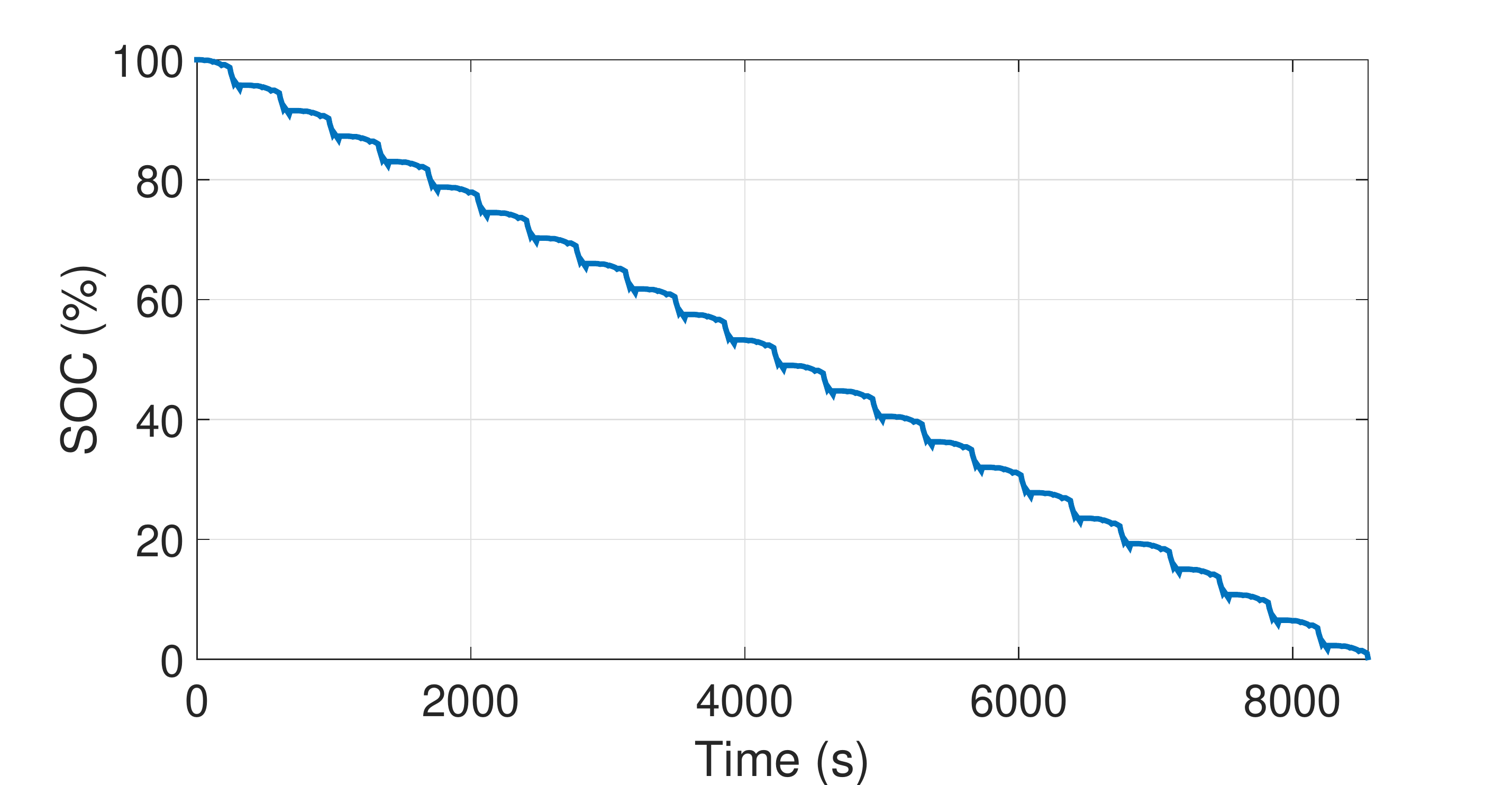}
\end{subfigure}
\caption{The open circuit voltage and SoC under the DST.}
\label{6}
\end{figure}

Fig. \ref{5} shows the fitting between the measured and estimated terminal voltage $V$ and the error. A slightly large deviation can be observed when the voltage is close to 2.0V, or when the SoC approaches zero. This is due to inaccuracies when measuring the open circuit voltage. The OCV-SoC lookup table and the estimated SoC are given in Fig. \ref{6}. The estimated SoC is used as a reference for the XKF.

\begin{figure}[h!]
  \centering
      \includegraphics[width=1\textwidth]{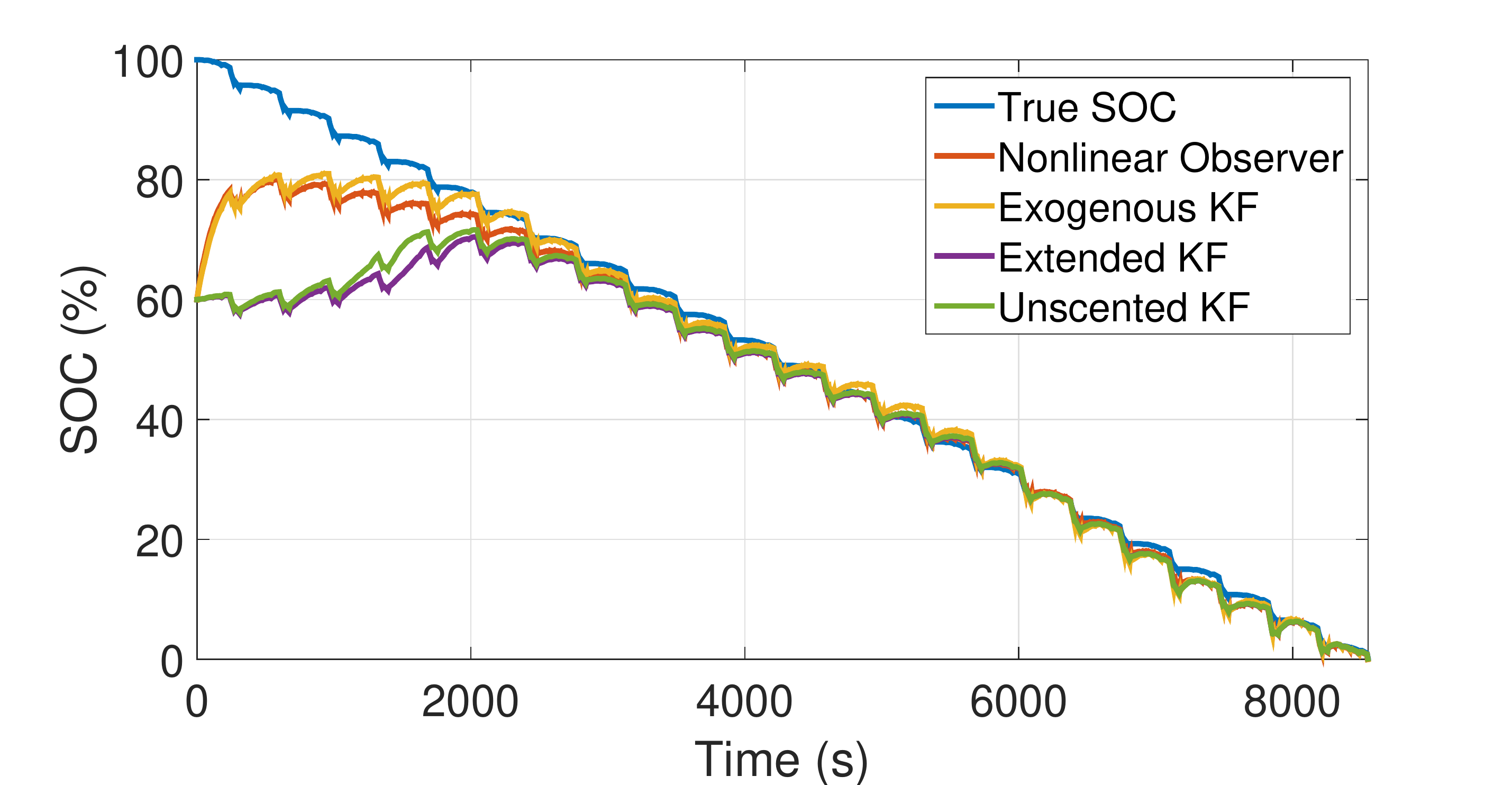}
  \caption{SoC estimation from different methods under DST with true initial condition 100\%.}
\label{7}
\end{figure}

\begin{figure}[h!]
  \centering
      \includegraphics[width=1\textwidth]{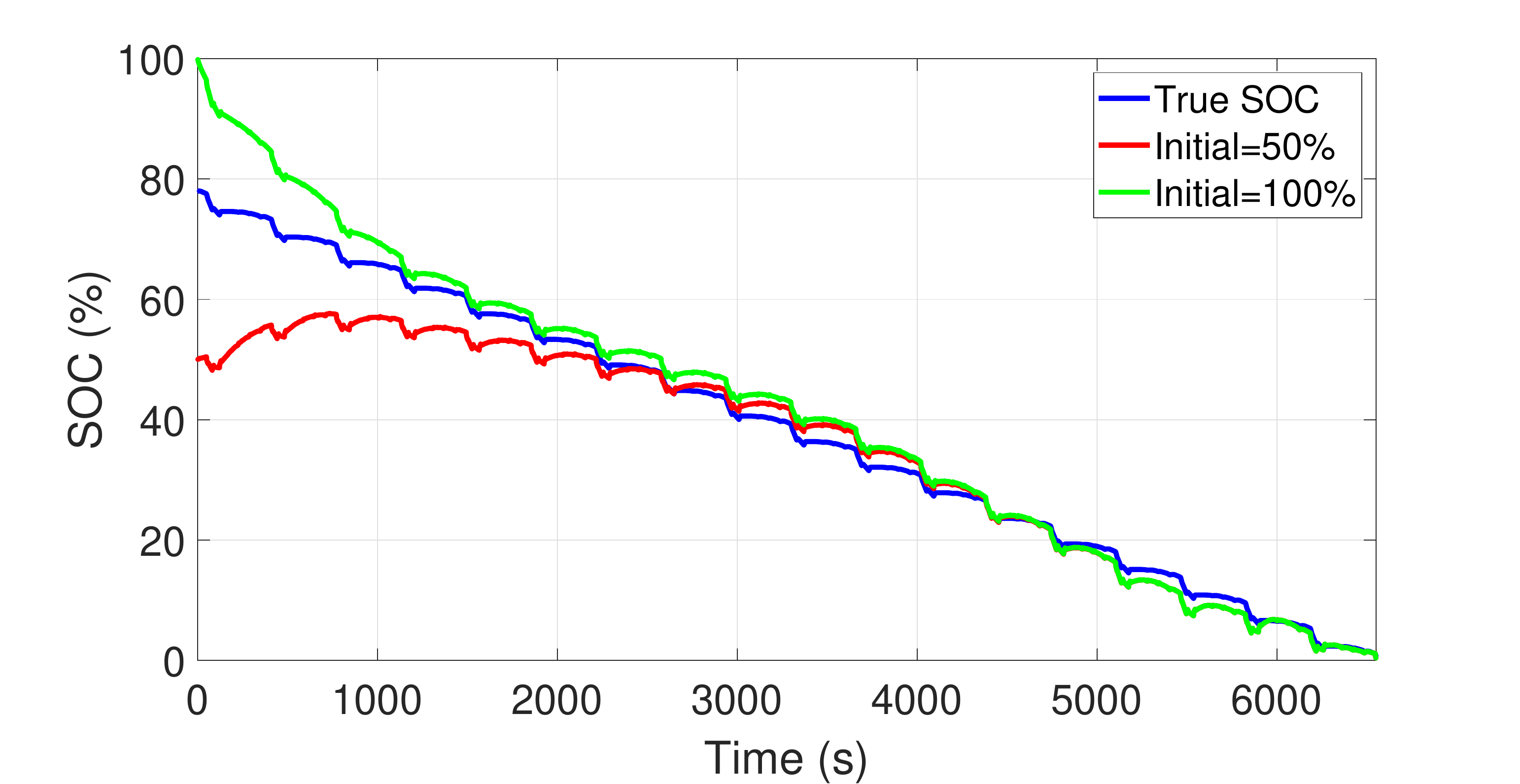}
  \caption{SoC estimation using XKF under DST with true initial condition 78\%.}
\label{71}
\end{figure}

To test the proposed method, we run the estimations with initial SoC at 60$\%$ and compare the results with some existing filter, e.g, EKF and UKF. The standard deviation for the process noise $\mathbf{Q}_n$ is 0.01, while for the measurement noise $\mathbf{R}_n$ is 0.04. These values are obtained using a simple covariance matching technique. The initial error covariance matrix $\mathbf{P}(0)$ is defined based on the initialization error. In this case, if the initial state is not very close, the value of $\mathbf{P}(0)$ should be large, whereas if the initialization is very good a smaller $\mathbf{P}(0)$ value can be used. The results are given in Fig. \ref{7}. It was observed that the NLO performs better in terms of the convergence rate than the EKF and UKF. The XKF further improves the estimation. Another experiment is done by starting the estimation when the SoC is at 78\%. Using XKF at different initial estimations, the values converge to the actual SoC, as can be seen from Fig. \ref{71}.

\begin{figure}[h!]
  \centering
      \includegraphics[width=1\textwidth]{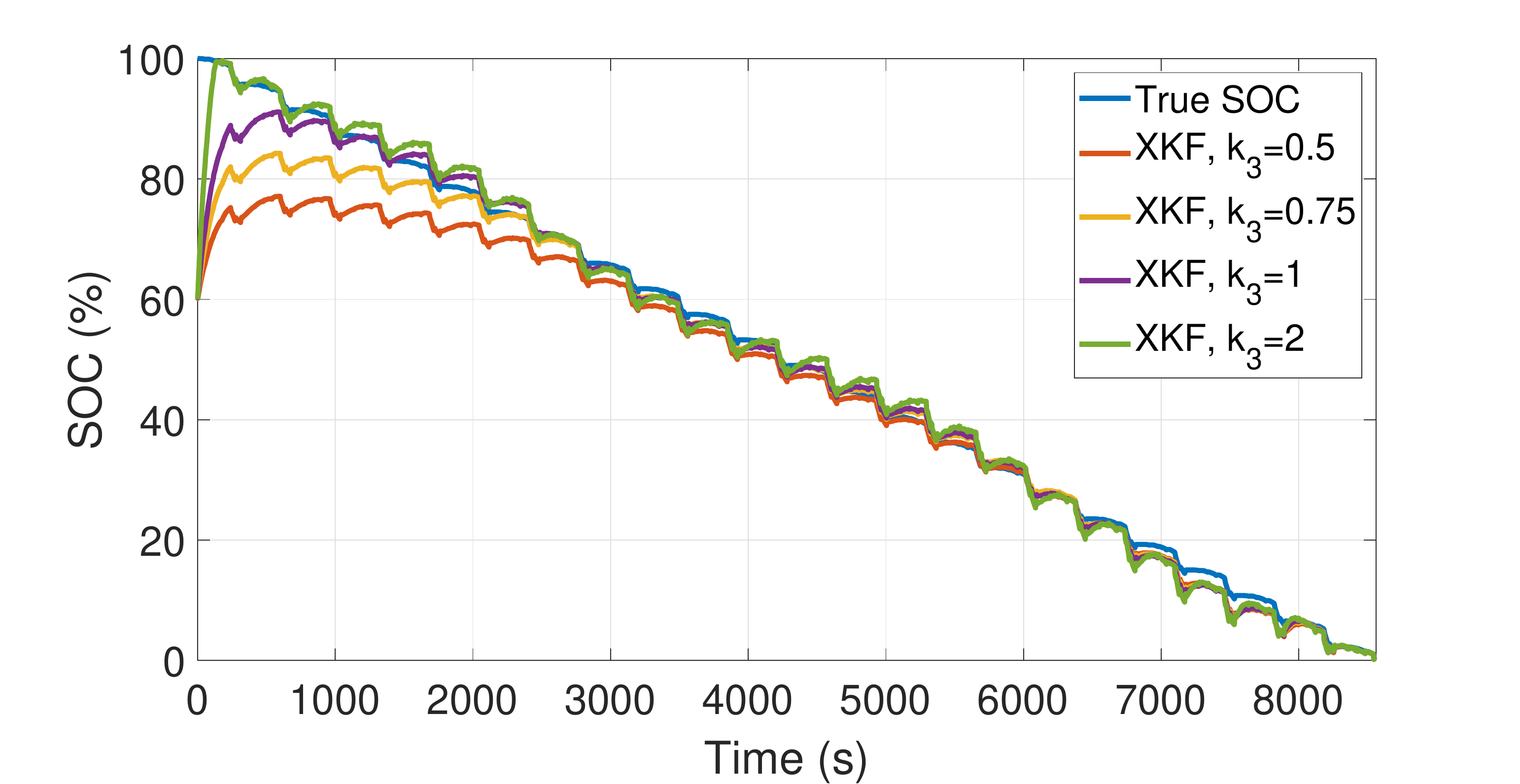}
  \caption{Comparison of different observer gain $k_3$.}
\label{hg}
\end{figure}

The convergence time can be increased using higher value of the gain $k_3$, as can be seen from Fig. \ref{hg}. With $k_3=2$ the estimate converges to the actual value in one minutes, while with $k_3=0.5$, the estimate converges after almost one hour. However, keep in mind that higher observer tends to overshoot as can be seen between $t=6500s$ to the end. Indeed, selecting an appropriate gain in crucial in SoC estimation. There is a trade-off between convergence rate and accuracy.

\subsection{Federal Urban Driving Schedule}

\begin{figure}[h!]
\centering
\begin{subfigure}[b]{1\textwidth}
   \includegraphics[width=1\linewidth]{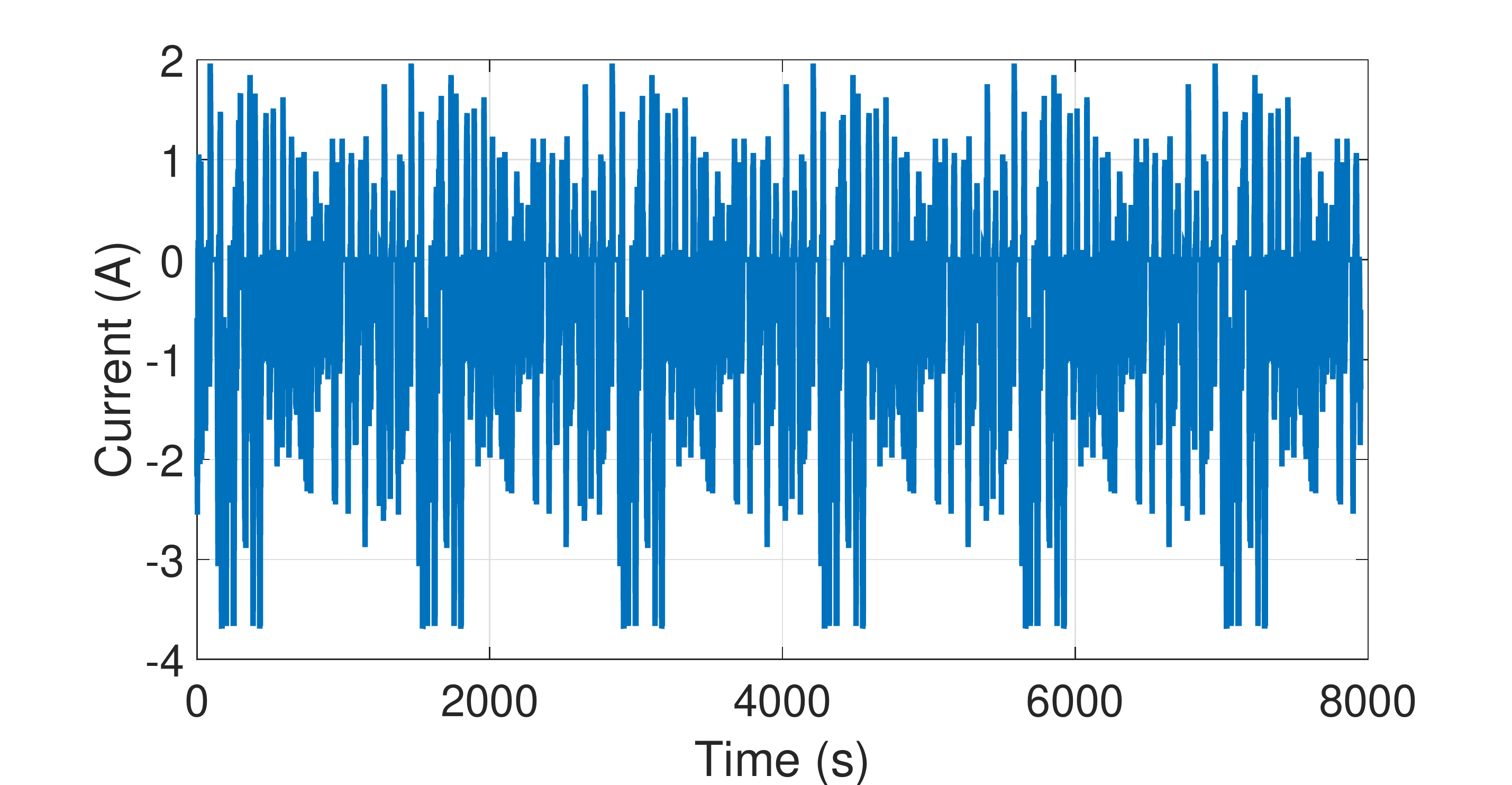}
\end{subfigure}
\begin{subfigure}[b]{1\textwidth}
   \includegraphics[width=1\linewidth]{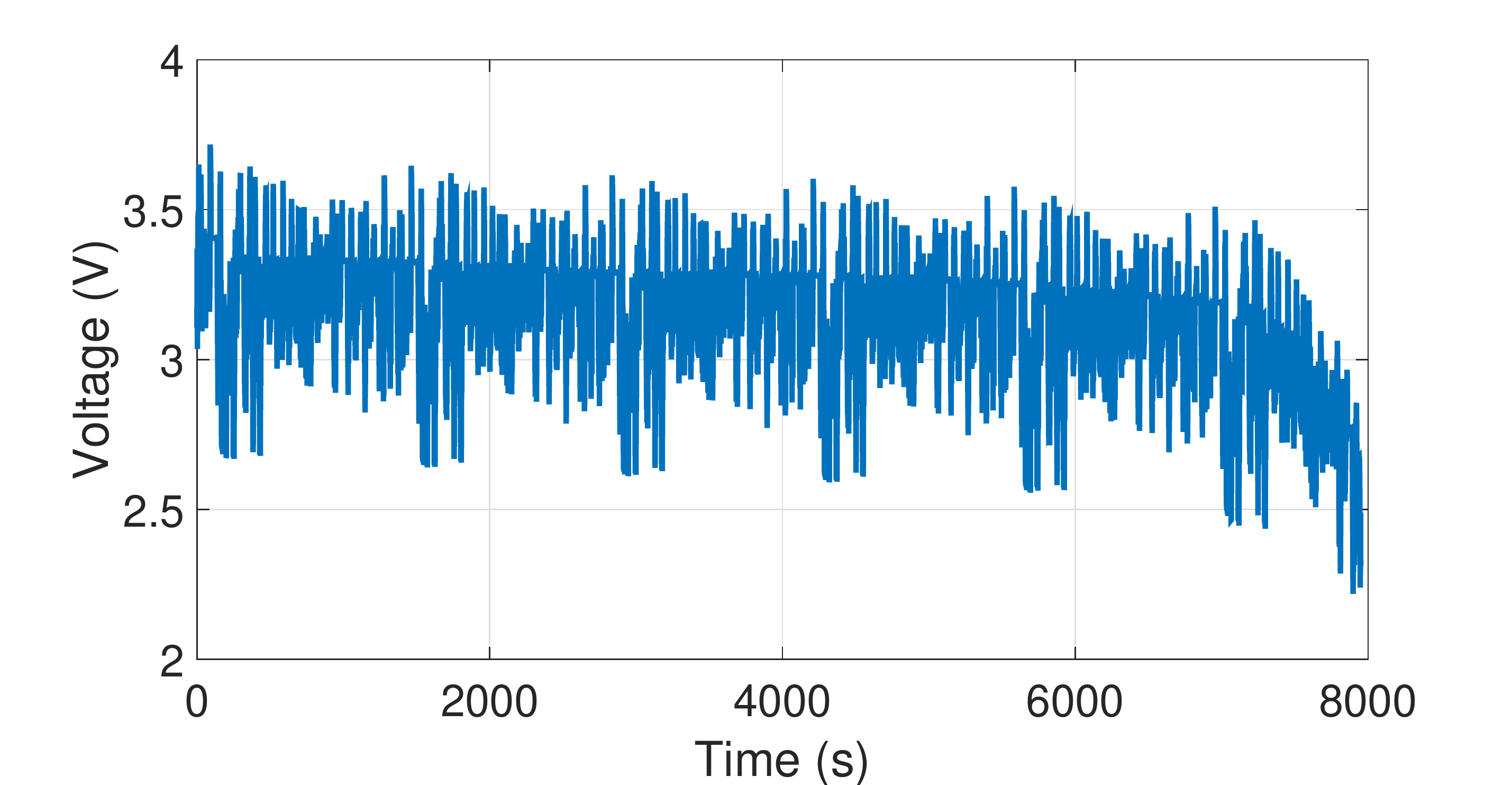}
\end{subfigure}
\caption{Measured current and terminal voltage response under the federal urban driving schedule (FUDS).}
\label{8}
\end{figure}

\begin{figure}[h!]
\centering
\begin{subfigure}[b]{1\textwidth}
   \includegraphics[width=1\linewidth]{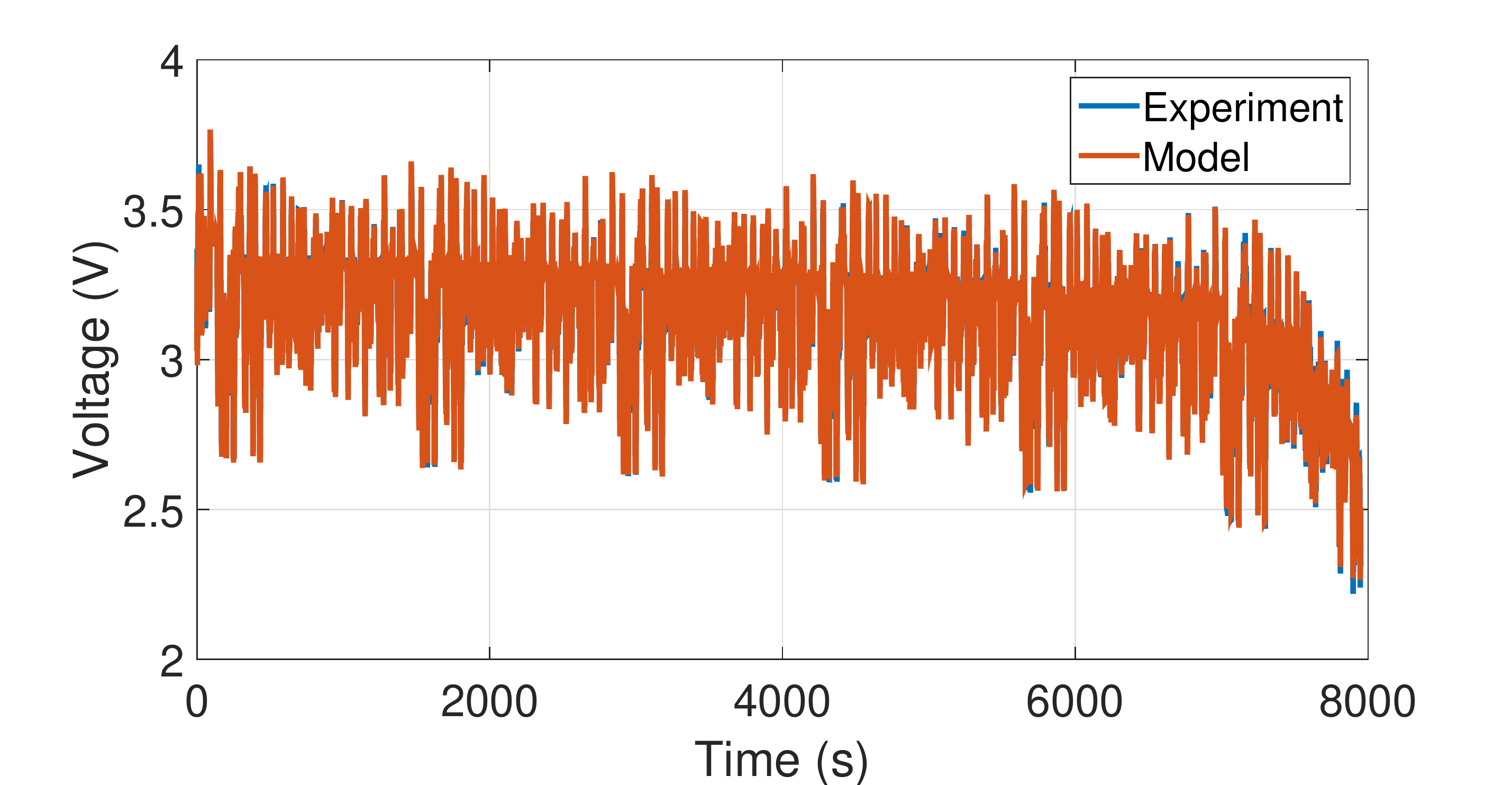}
\end{subfigure}
\begin{subfigure}[b]{1\textwidth}
   \includegraphics[width=1\linewidth]{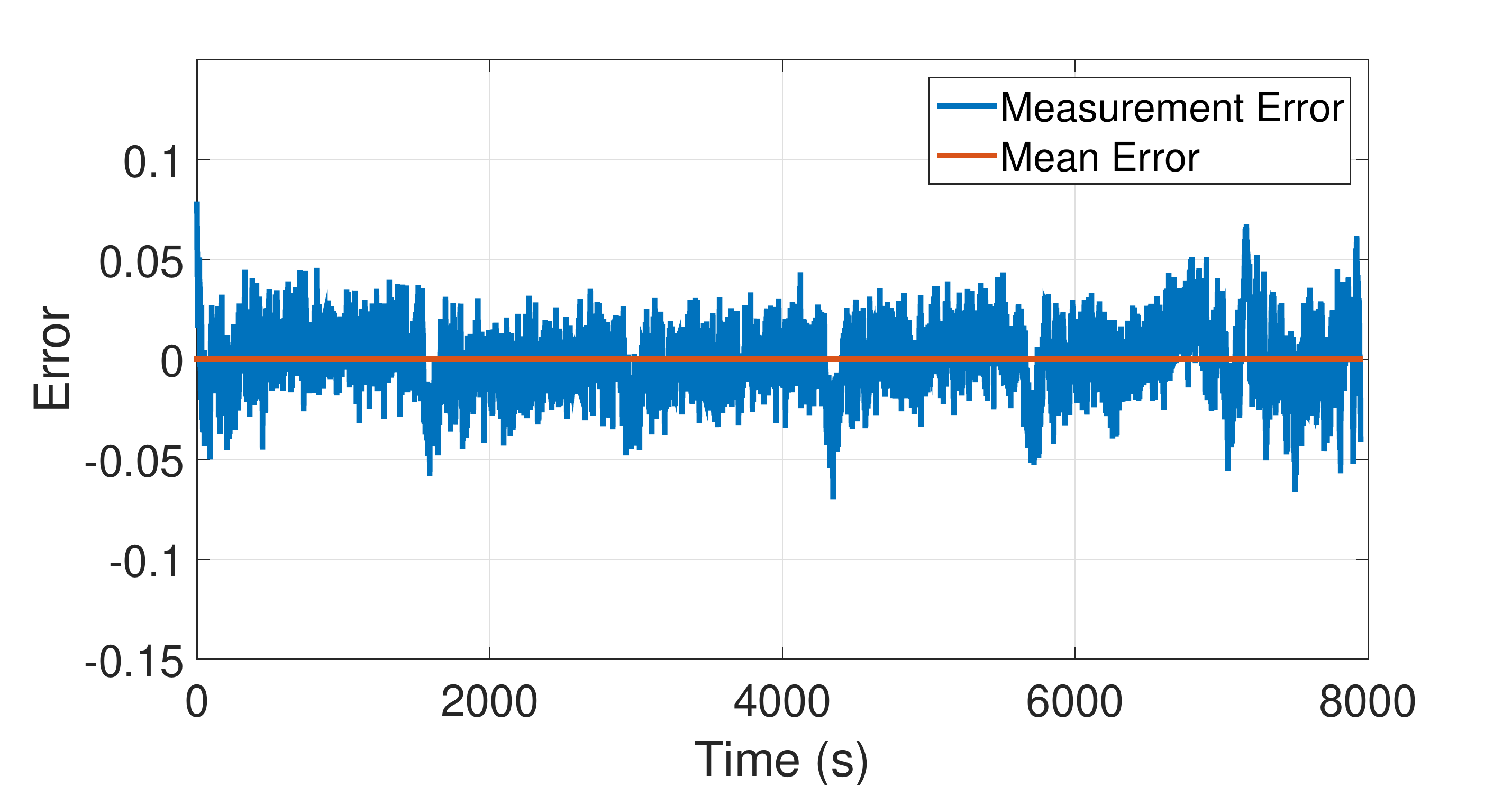}
\end{subfigure}
\caption{Measured and estimated voltage response under the FUDS and the model error.}
\label{10}
\end{figure}

\begin{figure}[h!]
\centering
\begin{subfigure}[b]{1\textwidth}
   \includegraphics[width=1\linewidth]{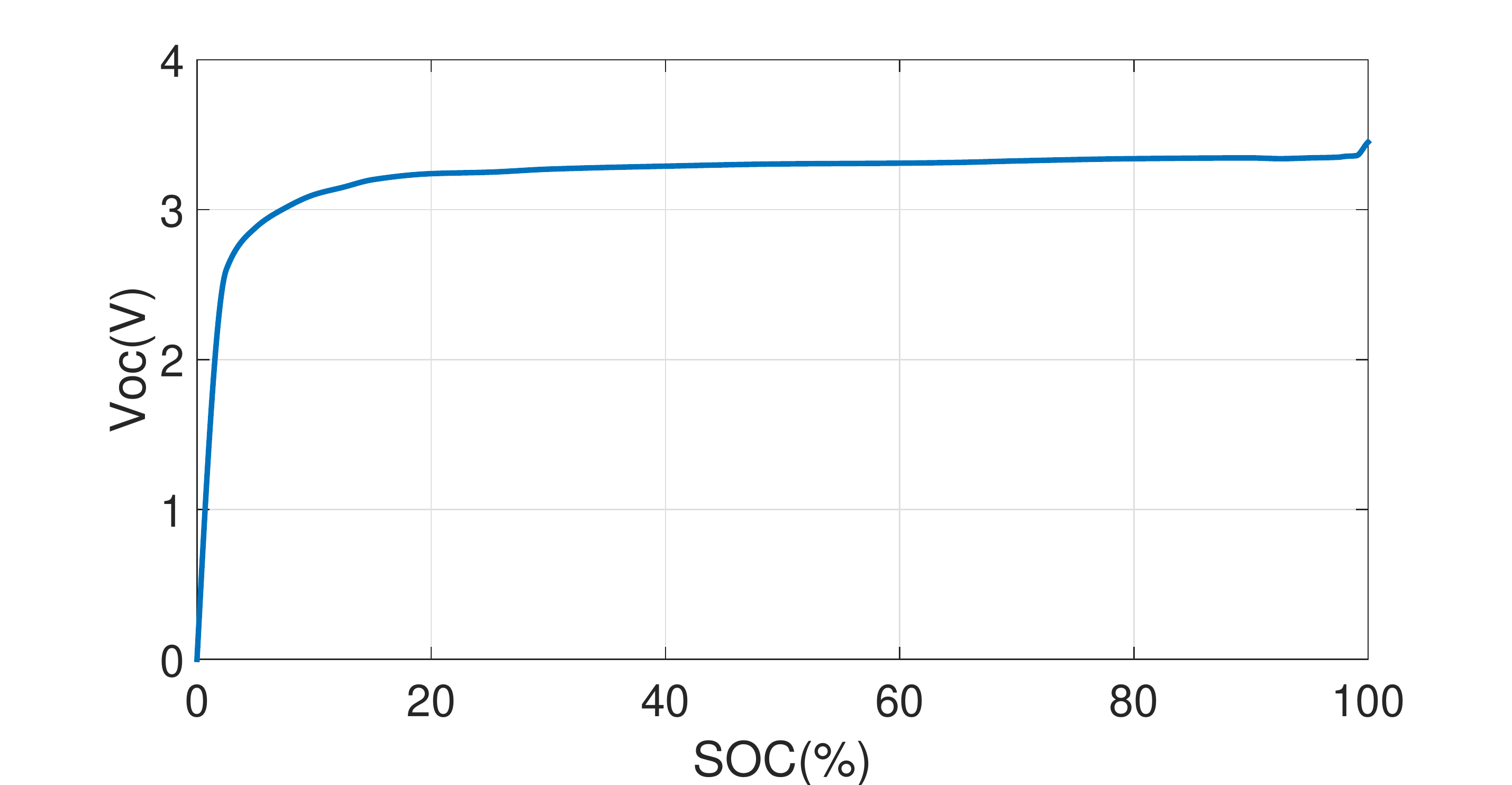}
\end{subfigure}
\begin{subfigure}[b]{1\textwidth}
   \includegraphics[width=1\linewidth]{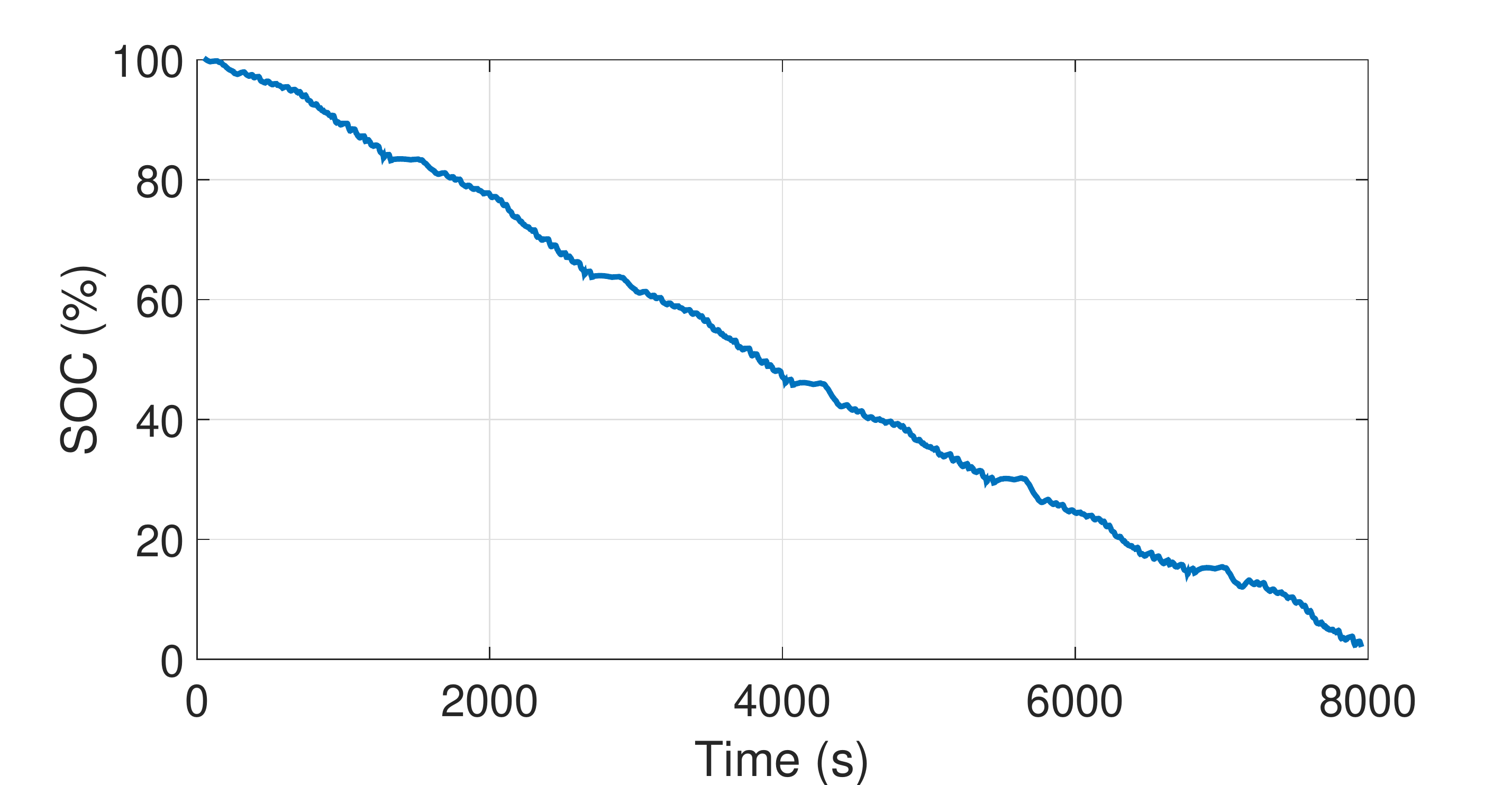}
\end{subfigure}
\caption{The open circuit voltage and SoC under the FUDS.}
\label{11}
\end{figure}

\begin{figure}[h!]
  \centering
      \includegraphics[width=1\textwidth]{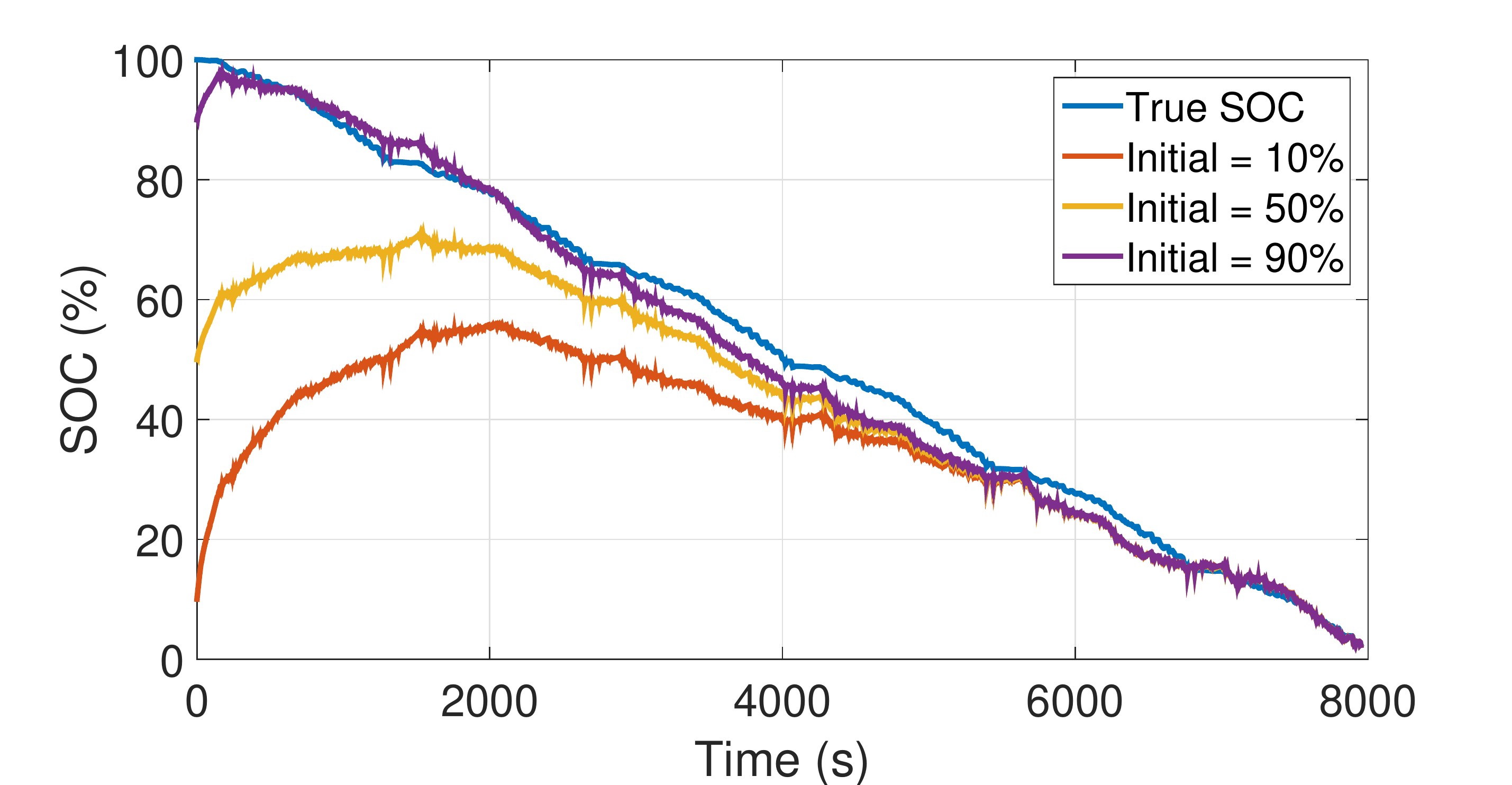}
  \caption{The SoC estimation from XKF under the FUDS at 20$^\circ$C.}
\label{9}
\end{figure}

A more realistic dynamic current test is given by the FUDS test. FUDS is based on a time-velocity profile from an automobile industry standard vehicle to test the dynamic electric vehicle performance. The measured current and voltage can be seen in Fig. \ref{8}. The model fitting and the error can be seen in Fig. \ref{10}, while the OCV-SoC lookup table and the SoC could be seen in Fig. \ref{11}, respectively. The current causes variation of the SoC from fully charged at 3.6V to empty at 2.0V. We run the simulations from three different initial guesses using XKF. The results are given in Fig. \ref{9}. It can be observed that the estimation converges to the SoC for any initial guess. The convergence rate can be increased using higher observer gain, but it  may cause significant overshoot and lead to oscillation and instability.

\subsection{Effect on Different Ambient Temperatures}

\begin{figure}[h!]
\centering
\begin{subfigure}[b]{1\textwidth}
   \includegraphics[width=1\linewidth]{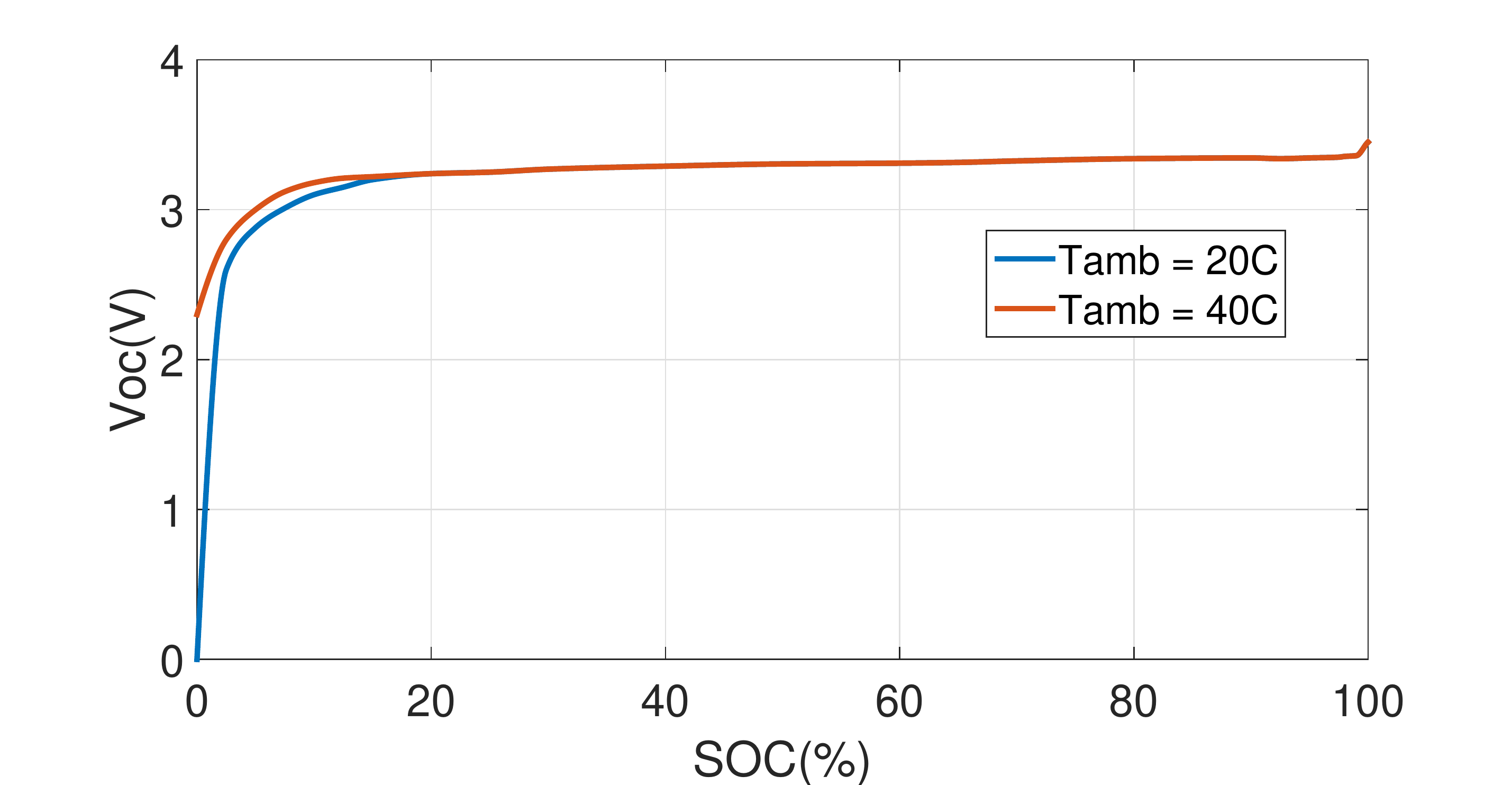}
\end{subfigure}
\begin{subfigure}[b]{1\textwidth}
   \includegraphics[width=1\linewidth]{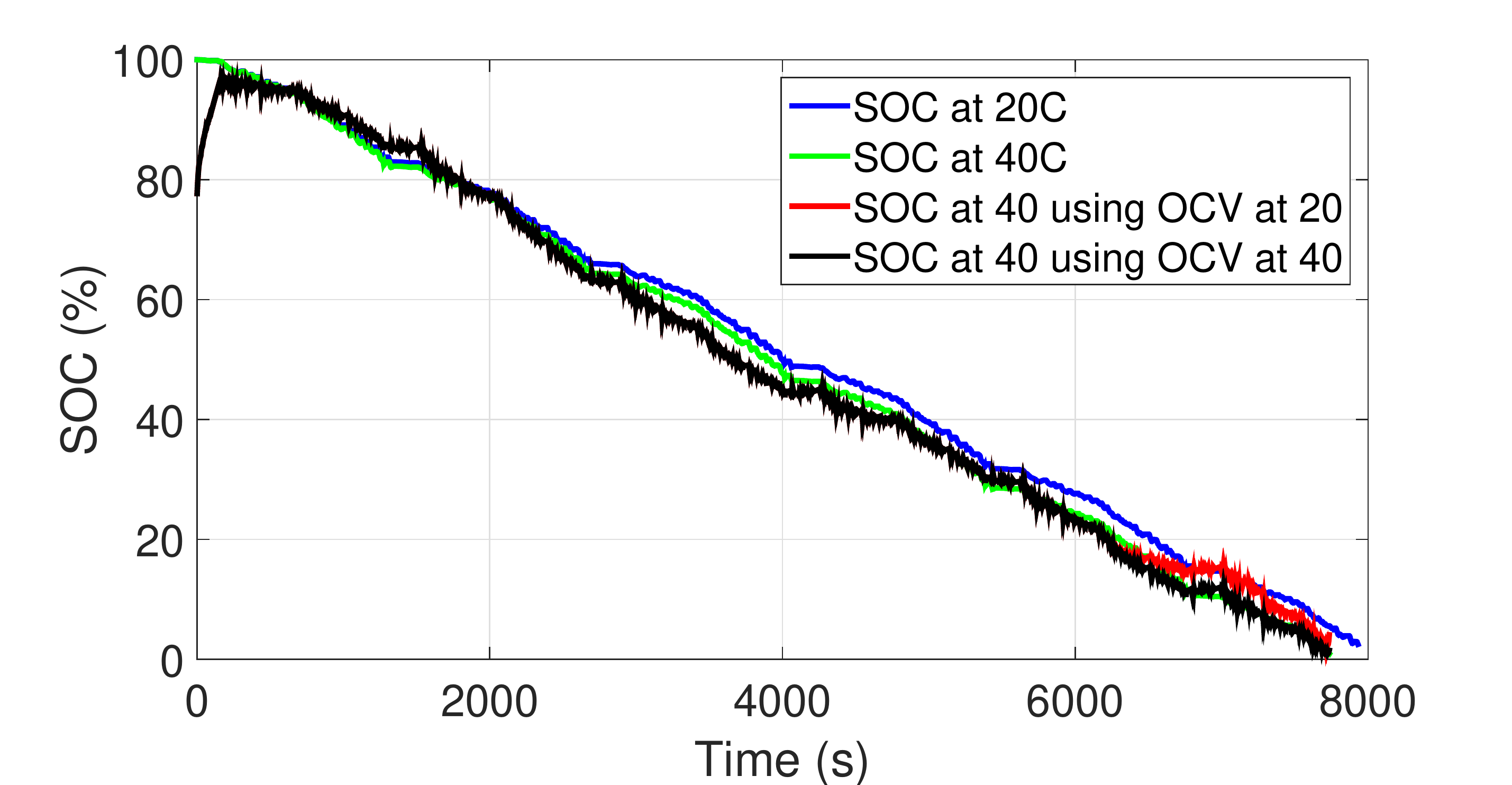}
\end{subfigure}
\caption{The open circuit voltage and SoC estimation using XKF under the FUDS at 40$^\circ$C.}
\label{14}
\end{figure}

In the last experiment, we estimate the SoC for different ambient temperatures. As we have discussed earlier, the OCV-SoC curve for LiFePO$_4$ has a wide flat region. Different ambient temperature will alter this curve significantly, especially at lower SoC. This can be seen from Fig. \ref{14}. The blue line shows the SoC at 20$^\circ$C, while the green line shows the SoC at 40$^\circ$C. We observe that higher ambient temperature will reduce the SoC. Furthermore, if the SoC estimation at 40$^\circ$C is calculated using the OCV-SoC curve at 20$^\circ$C (red line), then there is a deviation especially after the SoC is less than 20\%. While if it uses the OCV-SoC curve at 40$^\circ$C (black line), the estimate converges to the actual value. Thus, the OCV-SoC curve needs to be adjusted accordingly. Fail to update the OCV-SoC curve can cause incorrect SoC estimation. Thus, there is an incentive to model the temperature variation for the OCV-SoC curve, e.g, using a thermal-electrochemical model, which will be investigated in our further works.

\section{Conclusions and Future Works}

In this paper, we present SoC estimation of a LiFePO$_4$ battery using the eXogenous Kalman filter from the Coulomb counting and equivalent circuit model. The estimation algorithm consists of a cascade of NLO and LKF. Simulations against experimental data from the dynamics stress test and federal urban driving schedule show the algorithm is able to estimate the SoC accurately. Furthermore, comparisons with existing filter show the estimation using XKF converges faster, thanks to the exponential stability from the NLO. Future work includes the use of thermal-electrochemical model to improve the lithium battery estimation.

\section*{Acknowledgements}

This work was supported by the Free the Drones (FreeD) project at Center for Unmanned Aerial Vehicles, University of Southern Denmark, project number 5156-00008B, and the FRINATEK project “Multi-stage global sensor fusion for navigation” and the Centre of Autonomous Marine Operations and Systems (NTNU AMOS) at the Norwegian University of Science and Technology (NTNU), project numbers 250725 and 223254, respectively.

\end{document}